\documentclass[11pt,a4paper]{article}
\usepackage[T1]{fontenc}
\usepackage{lmodern,amsmath,amssymb,amsthm,mathtools}
\usepackage[margin=27mm]{geometry}
\usepackage{microtype,cite}
\usepackage[hidelinks]{hyperref}
\hypersetup{pdftitle={Global strong solutions of the three-dimensional quantum Navier--Stokes--Poisson system with large initial data},pdfauthor={Xiangdi Huang and Huitao Zhou}}
\allowdisplaybreaks[2]
\numberwithin{equation}{section}
\newtheorem{theorem}{Theorem}[section]
\newtheorem{proposition}[theorem]{Proposition}
\newtheorem{lemma}[theorem]{Lemma}
\theoremstyle{definition}
\newtheorem{definition}{Definition}[section]
\newcommand{\T}{\mathbb T^3}
\newcommand{\R}{\mathbb R}
\newcommand{\dd}{\,\mathrm d}
\newcommand{\budget}{\mathfrak E}
\DeclareMathOperator{\Div}{div}
\makeatletter
\renewcommand{\@fnsymbol}[1]{}
\makeatother
\title{Global strong solutions of the three-dimensional quantum Navier--Stokes--Poisson system with large initial data}
\author{
	\bfseries\large Xiangdi Huang\textsuperscript{a}\thanks{E-mail addresses: xdhuang@amss.ac.cn (X. Huang); zhouhuitao141@163.com (H. Zhou).}, Huitao Zhou\textsuperscript{a}\\
	\small a. Institute of Mathematics, Academy of Mathematics and Systems Science,\\
	\small Chinese Academy of Sciences, Beijing 100190, China;
}
\date{}
\begin{document}
\maketitle

\begin{abstract}
We prove that the three-dimensional quantum Navier--Stokes--Poisson system,
with density-proportional viscosity, matched viscous and quantum coefficients,
and pressure $P(\rho)=\rho^\gamma$, $1\le\gamma\le3$, has a unique global
strong solution with positive density on the unit torus $\T$ and in the whole
space $\R^3$, for arbitrarily large initial data. The constant background $b$
is the initial mean density on $\T$ and a prescribed positive far-field
density on $\R^3$. We assume $\inf\rho_0>0$, $r_0:=\rho_0-b\in H^3$ and
$u_0\in H^2$; on $\R^3$ we also require $E_0\in L^2$, where $E_0$ is the field
determined by Poisson's equation, equivalently $r_0\in\dot H^{-1}$. In both
geometries the density remains bounded above and below by positive constants
for all time, with constants depending only on the data and on the parameters
$\gamma$ and $\ell$. On $\T$, $(\rho,u,V)$ converges exponentially in
$H^3\times H^2\times H^5$ to the homogeneous equilibrium $(M,U,0)$ determined
by the conserved mass and momentum. The endpoint $\gamma=3$ is covered on $\T$
by a mixed density--velocity moment and on $\R^3$ by coupling the fourth
moment of the effective velocity $v=u+\nabla\log\rho$ to the fifth density
moment relative to $b$. Moments of $v$ of order $p\ge6$, together with a
parabolic comparison argument, yield the positive lower bound for the density.
In $\R^3$, the charge dissipation in the effective energy balance and
heat-kernel smoothing on short intervals give density bounds that are
uniform in time.
\end{abstract}



\section{Introduction}\label{sec:intro}

In a compressible quantum fluid with viscosity proportional to the density, the dissipation
produced by a given strain rate is small where the density is low. If the fluid is charged,
density fluctuations generate an electric field that acts back on the velocity. The question is
whether a large initial disturbance can force unbounded concentration of the density or drive it
towards vacuum. For the three-dimensional quantum Navier--Stokes--Poisson system with matched
viscous and quantum coefficients and pressure exponent $1\le\gamma\le3$,
Theorems~\ref{thm:periodic} and~\ref{thm:cauchy} give a unique global strong solution with
positive density for every datum in the classes specified in those theorems. In both geometries the density
remains bounded above and below by positive constants for all time, with constants depending only
on the data, and on the torus the solution converges exponentially to the homogeneous equilibrium
determined by its mass and momentum.

Let $\Omega$ be the unit torus $\T=\R^3/\mathbb Z^3$ or the whole space $\R^3$.
For a constant background density $b>0$ and an electric length $\ell>0$, the
density $\rho$, velocity $u$ and electrostatic potential $V$ solve
\begin{align}
 \rho_t+\Div(\rho u)&=0, \label{eq:mass}\\
 (\rho u)_t+\Div(\rho u\otimes u)+\nabla P(\rho)
 &=\Div(2\rho Du+\rho\nabla^2\log\rho)-\rho\nabla V, \label{eq:momentum}\\
 -\ell^2\Delta V&=\rho-b, \label{eq:poisson}
\end{align}
where $P(z)=z^\gamma$ with $1\le\gamma\le3$ is the barotropic pressure law.
The viscosity is proportional to the density, the bulk viscosity vanishes, and
the viscous and quantum coefficients are normalized to one.
We use $(a\otimes d)_{ij}=a_i d_j$, $(\nabla a)_{ij}=\partial_j a_i$,
$(\Div A)_i=\partial_j A_{ij}$ and $Da=(\nabla a+(\nabla a)^{\mathsf T})/2$,
with repeated spatial indices summed. Write $r=\rho-b$, $r_0=\rho_0-b$ and
$E=\nabla V$; the electric force in \eqref{eq:momentum} is then $-\rho E$.

On the torus, solvability of the Poisson equation fixes $b=M:=\int_{\T}\rho_0\dd x$;
we normalize $V$ to have zero spatial mean. In the whole space, $b$ is a
prescribed far-field density and we impose
\begin{equation}\label{eq:far}
   (\rho,u,V)(x,t)\to(b,0,0)\qquad(|x|\to\infty)
\end{equation}
for every fixed $t\ge0$.
These conventions fix the additive constant in $V$: the
zero spatial mean on $\T$, and the decay in \eqref{eq:far} on $\R^3$.

In one space dimension the large-data theory is comparatively complete, whereas
in higher dimensions the nonlocal electric field and the degenerate viscosity
add difficulties.
For the one-dimensional quantum Navier--Stokes system, J\"ungel obtained a global smooth solution
with strictly positive density and matched coefficients in the isothermal case \cite{J2011}. For unequal
viscosity and quantum coefficients, Chen and Zhao obtained global strong and classical solutions
for large positive-density data under conditions on the coefficients, the pressure exponent and the initial data \cite{CZ2025}. For the
one-dimensional quantum Navier--Stokes--Poisson system, Liu and Zhang proved global well-posedness
for large initial data and obtained upper and lower bounds, uniform in time, for the specific
volume $1/\rho$ and for the electric potential \cite{LiuZhang}.

In several dimensions, Danchin and Desjardins proved global existence for small perturbations of a
stable equilibrium for Korteweg fluids \cite{DD2001}. For the quantum Navier--Stokes--Poisson system
in $\R^3$, Wu and Hou proved the asymptotic stability of a stationary wave \cite{WuHou}, and Tong
and Xia obtained algebraic decay for small perturbations of a non-constant steady state with a small
doping profile \cite{TX2022}. For Navier--Stokes--Poisson--Korteweg flows, Xu and Li proved global
well-posedness for the system with density-dependent coefficients and small critical Besov data \cite{XuLi}. Gao and
Xu proved global existence of strong solutions for a class of large initial data for this system
\cite{GaoXu}. For the compressible Navier--Stokes--Poisson system without quantum stress, global
large-data results are known for spherically symmetric data on a ball \cite{FHL2026} and, in two
dimensions, for a general density-dependent bulk viscosity \cite{HP2026}.

For the density-dependent viscosity and quantum stress in \eqref{eq:momentum}, much of the
large-data strong theory has been developed for the system without the Poisson coupling. Haspot
stated a global existence result for the isothermal Cauchy problem \cite{Haspot2017}.
Huang, Meng and Zhang \cite[Remark~1.1]{HMZ2026} and Huang, Gu and Lei
\cite[Remark~1.1]{HGL2026} observed that Haspot's argument for a positive density lower
bound requires a weighted velocity estimate uniform in the moment order, and
that the estimates they obtain depend on that order. On the
three-dimensional torus, Huang, Meng and Zhang obtained global strong solutions for
$1\le\gamma<7/3$ \cite{HMZ2026}. Huang and Meng extended the range to $1\le\gamma<11/4$ for the
matched coefficients, obtaining a time-uniform upper bound for the density and a lower bound that
depends on the time horizon \cite{HM2026}. General density-dependent viscosity and capillarity
coefficients were treated in the periodic and whole-space settings \cite{GHMZ2026,HLZ2026}.

For the whole-space problem, the behavior of the density at infinity also enters the estimates. In
$\R^3$, Huang, Gu and Lei obtained global strong solutions for the quantum Navier--Stokes system in the range
$1\le\gamma<8/3$ with a positive far-field density \cite{HGL2026}. For initial data decaying at
infinity, Gu, Huang and Lei showed that on each finite time interval the density stays between
constant multiples of its initial spatial profile \cite{GHL2026}. In smooth bounded three-dimensional domains, Li, L\"u and
Zhou obtained global strong solutions for $1\le\gamma<3$ under a curvature Navier-slip condition on
the velocity and a homogeneous Neumann condition on the density \cite{LLZ2026}.

The results of this paper cover the full range $1\le\gamma\le3$, including the
cubic pressure $\gamma=3$, for the quantum Navier--Stokes--Poisson system with
a fixed positive background density. In both geometries, the upper and lower
density bounds are uniform in time and depend only on the data and on
$b,\gamma,\ell$. In $\R^3$, the extra requirement is the low-frequency
condition $r_0\in\dot H^{-1}$, equivalently $E_0\in L^2$; no smallness of the
data is needed. Relative to the results quoted above, the present results
include the endpoint $\gamma=3$ and give density bounds that are uniform in
time.

Norms and integrals without a domain are taken over $\Omega$. We write
$\|\cdot\|_p$ for the spatial $L^p$ norm, $H^k$ for the $L^2$ Sobolev space of
order $k$, $\dot H^s(\R^3)$ ($s\in\R$) for the homogeneous Sobolev space with
norm $\|f\|_{\dot H^s}=\||\xi|^s\widehat f\|_2$, and
$C_TX=C([0,T];X)$, $L_T^pX=L^p(0,T;X)$. The letter $C$ denotes a positive
constant whose value may change between estimates.

\begin{definition}[Positive strong solution]\label{def:positive-strong-solution}
A positive strong solution on $[0,T]$ is a distributional solution of
\eqref{eq:mass}--\eqref{eq:poisson} with initial traces $\rho(0)=\rho_0$ and
$u(0)=u_0$ and with the field $E=\nabla V$ determined by the Poisson equation,
such that $r\in C_TH^3$, $u\in C_TH^2$ and $E\in C_TL^2$, and such that the
continuous representative of $\rho$ is positive at every point of
$\Omega\times[0,T]$. On $\T$, we take $b=M$ and normalize $V$ to have zero
spatial mean. In the whole space we require \eqref{eq:far} for all $t\in[0,T]$.

On an interval $[a,a+h]$, the definition applies after translating the time
origin to $a$. A global positive strong solution is defined for all $t\ge0$
and satisfies the preceding conditions on every finite interval $[0,T]$.
\end{definition}

On the torus, compactness and pointwise positivity give a positive lower bound
for $\rho$ on $\T\times[0,T]$. In the whole space, the image of $t\mapsto r(t)$
is compact in $C_0(\R^3)$, the space of continuous functions vanishing at
infinity, so $\rho(x,t)\to b$ as $|x|\to\infty$ uniformly for $t\in[0,T]$.
Since $b>0$, pointwise positivity and continuity again give a positive lower
bound on $\R^3\times[0,T]$.

On the torus the total electric force vanishes,
$\int_{\T}\rho\nabla V\dd x=0$, and the momentum is conserved. The velocity of
the homogeneous equilibrium is therefore determined by the conserved
momentum, even for data with nonzero mean momentum.
\begin{theorem}\label{thm:periodic}
Let $\ell>0$, $1\le\gamma\le3$, and suppose that
\begin{equation}\label{eq:initial-data}
  \rho_0\in H^3(\T),\qquad u_0\in H^2(\T;\R^3),\qquad
  0<c_0\le\rho_0(x)\le C_0<\infty\quad(x\in\T).
\end{equation}
With $b=M$ and $\int_{\T} V\dd x=0$, the system
\eqref{eq:mass}--\eqref{eq:poisson} has a unique global positive strong
solution in the sense of Definition~\ref{def:positive-strong-solution}
with initial data $(\rho_0,u_0)$. For every finite $T>0$,
\begin{equation}\label{eq:strong-regularity}
\begin{aligned}
 \rho&\in C_TH^3\cap L_T^2H^4, & \rho_t&\in C_TH^1\cap L_T^2H^2,\\
 u&\in C_TH^2\cap L_T^2H^3, & u_t&\in C_TL^2\cap L_T^2H^1,\\
 V&\in C_TH^5\cap L_T^2H^6. &&
\end{aligned}
\end{equation}
There exists $C_{\rm den}\ge1$, depending only on the data, that is, on
$\ell$, $\gamma$, $c_0^{-1}$, $C_0$, $\|\rho_0\|_{H^3}$ and $\|u_0\|_{H^2}$,
such that $C_{\rm den}^{-1}\le\rho(x,t)\le C_{\rm den}$ for every
$(x,t)\in\T\times[0,\infty)$. Set $U=M^{-1}\int_{\T}\rho_0u_0\dd x$.
There are $T_0\ge0$, $C_{\rm dec}>0$ and $\sigma>0$ such that
\begin{equation}\label{eq:main-exponential-decay}
 \|\rho(t)-M\|_{H^3}+\|u(t)-U\|_{H^2}+\|V(t)\|_{H^5}
 \le C_{\rm dec} e^{-\sigma(t-T_0)},\qquad t\ge T_0.
\end{equation}
The time $T_0$ and the prefactor $C_{\rm dec}$ depend only on the data. The
rate $\sigma$ may be chosen in terms of $M,\gamma,\ell$ and the unit torus.
\end{theorem}

In the whole space, finite electric energy imposes a condition at low
frequencies. With the unitary Fourier transform
$\widehat f(\xi)=(2\pi)^{-3/2}\int_{\R^3}f(x)e^{-ix\cdot\xi}\dd x$,
Gauss law gives, for the field determined by the Poisson equation,
\begin{equation}\label{eq:data-field}
 \widehat E_0(\xi)=\ell^{-2}\frac{i\xi}{|\xi|^2}\widehat r_0(\xi),
 \qquad \ell^4\|E_0\|_2^2=
 \int_{\R^3}|\xi|^{-2}|\widehat r_0(\xi)|^2\dd \xi
 =\|r_0\|_{\dot H^{-1}}^2.
\end{equation}
Every $r_0\in L^1(\R^3)\cap H^3(\R^3)$ has finite initial electric energy:
$\widehat r_0$ is bounded near the origin, where $|\xi|^{-2}$ is locally integrable,
and $\widehat r_0\in L^2$ controls the integral over $|\xi|\ge1$.
In particular, integrable perturbations in this class may have nonzero total
charge $\int_{\R^3}r_0\dd x$.
\begin{theorem}\label{thm:cauchy}
Let $b,\ell>0$, $1\le\gamma\le3$, and suppose that
\begin{equation}\label{eq:data}
\begin{gathered}
  r_0\in H^3(\R^3)\cap\dot H^{-1}(\R^3),\qquad
  u_0\in H^2(\R^3;\R^3),\\
  0<c_0\le\rho_0(x)\le C_0<\infty\quad(x\in\R^3).
\end{gathered}
\end{equation}
Since $r_0\in H^3\subset C_0(\R^3)$, the bounds in \eqref{eq:data} force
$c_0\le b\le C_0$.
The system \eqref{eq:mass}--\eqref{eq:poisson} with far-field condition
\eqref{eq:far} has a unique global positive strong solution in the sense of
Definition~\ref{def:positive-strong-solution} with initial density $\rho_0$,
velocity $u_0$ and field $E_0$ determined by \eqref{eq:data-field}. For every
finite $T>0$,
\begin{equation}\label{eq:regularity}
\begin{aligned}
 r&\in C_TH^3\cap L_T^2H^4, & \rho_t&\in C_TH^1\cap L_T^2H^2,\\
 u&\in C_TH^2\cap L_T^2H^3, & u_t&\in C_TL^2\cap L_T^2H^1,\\
 E&\in C_TH^4\cap L_T^2H^5. &&
\end{aligned}
\end{equation}
There is a constant $C_{\rm den}\ge1$, depending only on
$b,\ell,\gamma,c_0^{-1},C_0,\|r_0\|_{H^3},\|u_0\|_{H^2}$ and
$\|r_0\|_{\dot H^{-1}}$, such that
\begin{equation}\label{eq:uniform}
   C_{\rm den}^{-1}\le\rho(x,t)\le C_{\rm den}
   \qquad ((x,t)\in\R^3\times[0,\infty)).
\end{equation}
\end{theorem}

The effective velocity $v=u+\nabla\log\rho$ is the key variable: it converts
the system into a parabolic one and makes the energy estimates possible. A
kinetic derivation of viscous corrections
to an isothermal quantum fluid model was given by Brull and M\'ehats \cite{BM2010}. Bresch,
Desjardins and Lin constructed large-data weak solutions of capillary flows with
density-proportional viscosity \cite{BreschDesjardinsLin}. The Bresch--Desjardins $\kappa$-entropy
yields an estimate for the density gradient \cite{BDZKappa}. Velocity integrability and
entropy estimates also underlie compactness results for degenerate barotropic fluids
\cite{MV2007,BVY2022}. In the quantum system, J\"ungel used the effective velocity
$u+\nabla\log\rho$ (in our normalization) to express the equations in parabolic form
\cite{J2010,J2011}. Augmented variables also yield relative-entropy estimates for comparing weak
and strong quantum flows \cite{BGL2019}.

The finite-energy weak theory allows vacuum, so the weak estimates give no positive lower bound for
the density. For periodic quantum Navier--Stokes flows, global existence of finite-energy weak
solutions was studied in \cite{AS2017,TZ2019}, and compactness in \cite{AS2018}. Weak solutions
have also been constructed with damping \cite{VY2016}, in a form admitting the semiclassical limit
\cite{LV2018}, and with non-trivial far-field behavior \cite{AHS2021}.

For Poisson-coupled quantum fluids, Yang and Ju studied weak solutions and the incompressible limit
\cite{YJ2015}, and Yang and Li obtained global weak solutions for the quantum
Navier--Stokes--Poisson system \cite{YL2017}. Quasineutral limits have been justified for viscous
capillary plasmas \cite{BDDPlasma,DM2015} and for quantum fluids with well-prepared data
\cite{LPW2017}. In the whole-space capillary model, dispersion controls oscillations of the
electric field \cite{DM2015}. With friction, weak solutions are available for large data on the
three-dimensional torus \cite{CC2025}, and the relaxation limit leads to quantum drift--diffusion
\cite{ACLS2021}. Since the viscosity is proportional to the density, a positive lower bound for the
density makes the viscosity nondegenerate; this is why the lower bound is central to the
continuation argument in Section~\ref{sec:continuation}.

To obtain the upper bound, we use the common diffusion in the equations for
$r$ and for the effective momentum $m=\rho v$, where $v=u+\nabla\log\rho$. In
particular, $\rho_t-\Delta\rho=-\Div m$. The continuity equation gives a
pointwise upper bound for $\rho$ once $m\in L^\infty(0,T;L^p)$ for some $p>3$.
The physical and effective energy balances provide the dissipative bounds
needed for this step: their dissipations control
$\int_0^T\|\nabla^2\sqrt\rho\|_2^2\dd t$ and hence the time integral of
$\|\rho\|_\infty^2$ on each finite window, independently of the unknown density
minimum (Lemmas~\ref{lem:quantum} and~\ref{lem:agmon} and
Proposition~\ref{prop:energy-bounds}).

On the torus, for $1<\gamma\le3$ we estimate density-weighted powers of the
effective specific energy $1+|v|^2/2+\rho^{\gamma-1}/(\gamma-1)$. The pressure
and electric work in Lemma~\ref{lem:mixed} contribute a growth coefficient
$C(1+\|\rho\|_\infty^{\gamma-1}+\|\rho\|_\infty)$, which is integrable in time
because $\gamma-1\le2$ and $\int_0^T\|\rho\|_\infty^2\dd t$ is bounded by
\eqref{eq:inputs}. At $\gamma=3$ the choice $a=3$ gives a sixth velocity moment
and a seventh density moment; together they yield
$m\in L^\infty(0,T;L^{7/2})$, so the time singularity of the differentiated
heat kernel is $h^{-13/14}$. At $\gamma=1$, separate velocity and density tests
give $m\in L^\infty(0,T;L^{16/5})$
(Proposition~\ref{prop:upper-iso}).

In the whole space, $\int_{\R^3}\rho^5\dd x=\infty$ because the background density is positive. We
pair the fourth velocity moment $\int\rho|v|^4\dd x$ with the relative fifth density moment
\[
   \int_{\R^3}\bigl[\rho^5-b^5-5b^4(\rho-b)\bigr]\dd x.
\]
At $\gamma=3$, both the fourth-moment estimate and the fifth-moment estimate involve
$\int\rho^5|v|^2\dd x$. On $\{\rho>\max(1,2b)\}$ we have
$\rho^5\le(16/13)\Phi(\rho)$, where $\Phi(\rho)=\rho^5-b^5-5b^4(\rho-b)$. Integrating the identity
\eqref{eq:joint-square} over this set therefore bounds
$\int_{\{\rho>\max(1,2b)\}}\rho^5|v|^2\dd x$ by $C\|\rho\|_\infty^2$ times the sum of the two
moments. On a window of prescribed length, Proposition~\ref{prop:upper} thus requires only upper
bounds on $\rho$ and $|v|$ at the initial time, together with the original energy bound and the
fixed parameters.

Ruling out vacuum requires a different use of the velocity moments. The
estimates of \cite{HMZ2026} depend logarithmically on a lower bound for the
density, which suggests retaining the dependence on the moment order. With
finite-interval bounds on $\|\rho\|_\infty$ and $\|E\|_\infty$,
Lemma~\ref{lem:large-moments} bounds the weighted moment
$(\int\rho|v|^p\dd x)^{1/p}$ by $A_T\sqrt p$ for every real $p\ge6$. On each
short interval we choose the moment order $p=6+\log\kappa$, where
$\kappa=\max\{e,\sup\|\rho^{-1}\|_\infty\}$; removing the density weight then
costs only the factor $\kappa^{1/p}\le e$. On subintervals of length
proportional to $1/p$, a parabolic
comparison argument preserves a fixed fraction of the density minimum at the
start of the subinterval. Iterating this estimate yields a bound for the
reciprocal density on each finite interval (Section~\ref{sec:lower}). An eighth
moment of the effective velocity bounds the convective flux in the effective
momentum equation. Heat smoothing bounds $m$ pointwise, and the resulting
bounds on $m$ and $\rho$ permit continuation of the Sobolev solution in
Section~\ref{sec:continuation}.

It remains to make the density bounds independent of the time horizon. On the torus, relative
entropy has been used to prove exponential decay for Korteweg flows with drag \cite{BGLV2022}; here
we use instead the dissipation and the conserved momentum, which yield a time at which
$(\rho-M,v-U)$ is small in $H^2\times H^1$. After the change of variables
that moves with $U$, the spectral gap of the
linearized operator and a coupled estimate for the pressure and electric terms give
\eqref{eq:main-exponential-decay} (Section~\ref{sec:relaxation}).

In the whole space, the time integrals of $\|r\|_2^2$,
$\|\nabla^2\sqrt\rho\|_2^2$ and $\|\sqrt\rho\,\nabla v\|_2^2$ are bounded by
\eqref{eq:dissipation-budget}. These bounds give a fixed window length
$T_{\rm w}$ such that, for each observation time $t\ge T_{\rm w}$, the argument
in Section~\ref{sec:main-proof} yields a $t_0\in[t-T_{\rm w},t-1]$ at which
$9b/16\le\rho(x,t_0)\le25b/16$ for all $x\in\R^3$ and $\|m(t_0)\|_6\le1$.
After a fixed delay $\delta>0$, heat smoothing gives an $L^\infty$ bound for
$v(t_0+\delta)$ independent of the position of the interval. From this time on,
the upper and
lower density bounds propagate to $t$ with constants independent of $t$.
Appendix~\ref{app:local} proves local existence and uniqueness in the class
of Definition~\ref{def:positive-strong-solution}.

\section{Effective equations and energy estimates}\label{sec:structure}
Set $s=\sqrt\rho$, $\Pi=P(\rho)-P(b)$, and
\begin{equation}\label{eq:effective-variables}
   v=u+\nabla\log\rho,\qquad m=\rho v=\rho u+\nabla r,\qquad
   \mathcal U=(r,m,E).
\end{equation}
Both $r$ and $m$ solve heat equations with the same principal part
$-\Delta$. The field can be included as a third heat unknown provided the
initial Poisson constraints are imposed.

\subsection{Effective momentum and reconstruction of the field}\label{sec:effective}
Adding the gradient of \eqref{eq:mass} to \eqref{eq:momentum} gives
\begin{equation}\label{eq:effective-flux}
 \partial_t m_i+\partial_i\Pi+\rho E_i
   =\partial_j(\rho\,\partial_jv_i-m_i u_j).
\end{equation}
Indeed, in components,
\begin{equation}\label{eq:effective-cancellation}
\begin{aligned}
\bigl(\Div(2\rho Du)-\nabla\Div(\rho u)\bigr)_i
&=\partial_j\bigl(\rho\partial_ju_i-u_j\partial_i\rho\bigr),\\
\bigl(\Div(\rho\nabla^2\log\rho)-\Div(\rho u\otimes u)\bigr)_i
&=\partial_j\bigl(\rho\partial_i\partial_j\log\rho-\rho u_iu_j\bigr),
\end{aligned}
\end{equation}
where we used $\Div(2\rho Du)_i=\partial_j(\rho\partial_ju_i+\rho\partial_iu_j)$
and $\partial_i(\rho u_j)=\rho\partial_iu_j+u_j\partial_i\rho$. Adding the two
identities and using $\partial_jv_i=\partial_ju_i+\partial_i\partial_j\log\rho$,
$\partial_i\rho=\partial_ir$ and $m_i=\rho u_i+\partial_ir$ gives
\eqref{eq:effective-flux}.
Since $\rho\,\partial_jv_i-m_i u_j=\partial_j m_i-m_i v_j$, the resulting system is
\begin{align}
 r_t-\Delta r&=-\Div m,\label{eq:heat-mass}\\
 m_t-\Delta m&=-\Div(m\otimes m/\rho)-\nabla\Pi-\rho E. \label{eq:heat-momentum}
\end{align}
The continuity equation also gives
\begin{equation}\label{eq:effective-velocity}
 \rho(v_t+u\cdot\nabla v)=\Div(\rho\nabla v)-\nabla\Pi-\rho E.
\end{equation}
Dividing by $\rho$ and using
$\Div(\rho\nabla v)/\rho=\Delta v+\nabla\log\rho\cdot\nabla v$ and
$u=v-\nabla\log\rho$, we obtain
\begin{equation}\label{eq:v}
 v_t-\Delta v+(v-2\nabla\log\rho)\cdot\nabla v+
       \gamma\rho^{\gamma-2}\nabla\rho=-E. 
\end{equation}

Let $Q$ denote the orthogonal projection in $L^2(\Omega;\R^3)$ onto the
closure of gradients, with symbol $\xi\otimes\xi/|\xi|^2$ for $\xi\ne0$.
On $\T$ the frequencies are $2\pi\mathbb Z^3$, the Fourier coefficients
are $\widehat f(\xi)=\int_{\T}f(x)e^{-i x\cdot\xi}\dd x$, and the zero mode
of $Q$ is zero. Gauss law and $\nabla\times E=0$ give
\begin{equation}\label{eq:gauss-symbol}
 \widehat E(\xi)=\ell^{-2}\frac{i\xi}{|\xi|^2}\widehat r(\xi)
 \quad(\xi\ne0),\qquad QE=E,\qquad \|\nabla E\|_2=\ell^{-2}\|r\|_2.
\end{equation}
The Fourier identity is understood almost everywhere in $\R^3$.
In the periodic problem $r$ and $E$ have zero mean.
The field $E_0+\ell^{-2}\int_0^t Q(\rho u)(\tau)\dd\tau$ lies in the range of
$Q$, and it has the same divergence and curl as $E(t)$: its divergence is
\[
-\ell^{-2}r_0+\ell^{-2}\int_0^t\Div(\rho u)(\tau)\dd\tau=-\ell^{-2}r(t)
\]
by \eqref{eq:mass}, and its curl vanishes. Uniqueness in the range of $Q$ gives
$E_t=\ell^{-2}Q(\rho u)$. Since
$\Delta E=-\ell^{-2}\nabla r$ and $Q\nabla r=\nabla r$, this is equivalent to
\begin{equation}\label{eq:field-heat}
   E_t-\Delta E=\ell^{-2}Qm.
\end{equation}
Conversely, the system \eqref{eq:heat-mass}--\eqref{eq:heat-momentum},
\eqref{eq:field-heat} propagates the constraints
\begin{equation}\label{eq:constraint-heat}
 (\partial_t-\Delta)(\Div E+\ell^{-2}r)=0,\qquad
 (\partial_t-\Delta)(\nabla\times E)=0.
\end{equation}
For $\mathcal U\in C_TH^2$ these expressions in parentheses belong to
$C_TH^1$. Zero initial data give zero solutions. In the periodic setting
\eqref{eq:heat-mass} and \eqref{eq:field-heat} also preserve the zero means.
For positive $\rho$ and $r\in C_TH^3$ we recover the physical equations
by reversing \eqref{eq:effective-flux} with $u=(m-\nabla r)/\rho$.

On the torus, the potential is given by the mean-zero inverse Laplacian.
In $\R^3$, truncate the multiplier $-i\xi\cdot\widehat E/|\xi|^2$ to
$1/n\le|\xi|\le n$. The gradients of the resulting potentials converge to
$QE=E$ in $L^2$, so the homogeneous Sobolev inequality gives convergence of
the potentials in $L^6$. Since $E\in H^2\subset L^\infty$, the limit $V$ is
Lipschitz continuous. If $V$ did not tend to zero at infinity, its Lipschitz
continuity would produce infinitely many disjoint balls on which $|V|$ is
bounded below, contradicting $V\in L^6$. Thus the whole-space normalization
takes the form
\begin{equation}\label{eq:potential-normalized}
  \nabla V=E,\qquad V\in L^6(\R^3),\qquad \lim_{|x|\to\infty}V(x,t)=0.
\end{equation}

For later estimates write $\mathcal U_t-\Delta\mathcal U
=-\Div\mathcal F(\mathcal U)+\mathcal S(\mathcal U)$ with
\begin{equation}\label{eq:heat-maps}
 \mathcal F(\mathcal U)=(m,m\otimes m/\rho+\Pi I_3,0),\qquad
 \mathcal S(\mathcal U)=(0,-\rho E,\ell^{-2}Qm).
\end{equation}
Here $I_3$ is the identity matrix. Let $G_h(x)=(4\pi h)^{-3/2}e^{-|x|^2/(4h)}$
in $\R^3$ and use its periodization on $\T$. On either domain the heat
equation $\zeta_t-\Delta\zeta=-\Div F+g$ has the representation
\begin{equation}\label{eq:heat-integral}
 \zeta(t)=G_{t-a}*\zeta(a)+\int_a^t\!
 \left(G_{t-\tau}*g(\tau)-\sum_j\partial_jG_{t-\tau}*F_{\cdot j}(\tau)\right)d\tau.
\end{equation}
In $\R^3$, fix the constant
\begin{equation}\label{eq:gaussian-constant}
 C_{\rm heat}=8\left(1+\sup_{1\le q\le\infty}
 \left[\|G_1\|_q+\sum_{j=1}^3\|\partial_jG_1\|_q\right]\right).
\end{equation}
For $1\le p\le q\le\infty$, $1+1/q=1/k+1/p$ and $j=0,1$, scaling gives
\begin{equation}\label{eq:heat-kernel}
 \|\nabla^jG_h*g\|_q\le h^{-j/2-(3/2)(1/p-1/q)}
       \|\nabla^jG_1\|_k\|g\|_p \qquad(\Omega=\R^3).
\end{equation}
Each component is estimated separately when a derivative is involved.
On the unit torus the same powers of $h$ hold for $0<h\le1$ with an
absolute constant independent of $p,q$. Indeed the periodized kernel and
its derivatives have $L^1$ bounds $C h^{-j/2}$ by summation, and
$L^\infty$ bounds $C h^{-(3+j)/2}$ by Gaussian decay. Interpolation and
Young's inequality give the assertion.

For $\zeta(a)\in H^k$, $F,g\in L^2(a,T;H^k)$ and $k\ge0$, the $H^k$ inner product
with $\zeta$ and Young's inequality give
\begin{equation}\label{eq:linear-energy}
 \frac{d}{dt}\|\zeta\|_{H^k}^2+\|\nabla\zeta\|_{H^k}^2
 \le\|F\|_{H^k}^2+\|g\|_{H^k}^2+\|\zeta\|_{H^k}^2.
\end{equation}
In particular,
\begin{equation}\label{eq:linear-integrated}
 e^{-(t-a)}\|\zeta(t)\|_{H^k}^2+
 \int_a^t e^{-(\tau-a)}\|\nabla\zeta\|_{H^k}^2\dd\tau
 \le\|\zeta(a)\|_{H^k}^2+
 \int_a^t e^{-(\tau-a)}(\|F\|_{H^k}^2+\|g\|_{H^k}^2)\dd\tau.
\end{equation}
Fourier truncation proves this estimate for the mild solution and yields
$\zeta\in C([a,T];H^k)\cap L^2(a,T;H^{k+1})$.
We use $J=(1-\Delta)^{1/2}$, $\Lambda=(-\Delta)^{1/2}$ and
$\|w\|_{H^k}=\|J^kw\|_2$. An $L^2$ inner product is spatial and acts on
all components. With $L^2$ as pivot space we denote $H^{-1}$--$H^1$
duality by $\langle\,\cdot\,,\,\cdot\,\rangle_{H^{-1},H^1}$.
In particular $w\in L^2H^{k+1}$ and $w_t\in L^2H^{k-1}$ imply
$w\in C H^k$ and $\tfrac12(d/dt)\|w\|_{H^k}^2
=\langle J^kw_t,J^kw\rangle_{H^{-1},H^1}$ almost everywhere, by
Fourier truncation. Time norms here and below refer to the indicated
interval. We write $A\asymp B$ for two-sided bounds by positive multiples
with the stated parameter dependence.

\subsection{Energy balances and the quantum Hessian}\label{sec:energy}
For $z\ge0$ let the relative internal energy be
\begin{equation}\label{eq:relative-functions}
 H(z\mid b)=
 \begin{cases}
 [z^\gamma-b^\gamma-\gamma b^{\gamma-1}(z-b)]/(\gamma-1),&\gamma>1,\\
 z\log(z/b)-z+b,&\gamma=1.
 \end{cases}
\end{equation}
Use its continuous value at zero and abbreviate it to $H(z)$ at fixed $b$.
Then $H(b)=H'(b)=0$ and $zH''(z)=P'(z)$ for $z>0$.
In the definitions of $E_u$ and $E_v$, take $U=M^{-1}\int_\T\rho_0u_0\dd x$
on the torus and $U=0$ in the whole space. Define
\begin{align}
 E_u(t)&=\int_\Omega\left(\tfrac12\rho|u-U|^2+2|\nabla s|^2+H(\rho)
                    +\tfrac{\ell^2}{2}|E|^2\right)\dd x, \label{eq:energy-u}\\
 E_v(t)&=\int_\Omega\left(\tfrac12\rho|v-U|^2+H(\rho)
                    +\tfrac{\ell^2}{2}|E|^2\right)\dd x . \label{eq:energy-v}
\end{align}
The physical energy dissipation controls $\sqrt\rho\,Du$, and the effective
energy dissipation controls $\sqrt\rho\,\nabla v$, both in $L^2(0,T;L^2)$.
Since $Dv-Du=\nabla^2\log\rho$, the two estimates together bound
$\sqrt\rho\,\nabla^2\log\rho$ in the same space. Lemma~\ref{lem:quantum}
then bounds $\nabla^2\sqrt\rho$ in $L^2(0,T;L^2)$.

\begin{lemma}\label{lem:energy}
On each compact interval of positive strong existence $E_u$ and $E_v$ are
absolutely continuous and satisfy, almost everywhere,
\begin{align}
 E_u'+2\int_\Omega\rho|Du|^2\dd x&=0,\label{eq:Eu}\\
 E_v'+\int_\Omega\rho|\nabla v|^2\dd x+
       \gamma\int_\Omega\rho^{\gamma-2}|\nabla\rho|^2\dd x
       +\ell^{-2}\|r\|_2^2&=0 . \label{eq:Ev}
\end{align}
On the torus $\int_\T\rho\dd x=M$ and $\int_\T\rho u\dd x=\int_\T\rho v\dd x=MU$.
\end{lemma}
\begin{proof}
Integrating \eqref{eq:mass} over $\T$ gives mass conservation. Since
$\rho=b-\ell^2\Delta V$, the integral of the electric force vanishes:
\[
\int_\T\rho\,\partial_iV\dd x
=-\ell^2\int_\T\Delta V\,\partial_iV\dd x
=+\frac{\ell^2}{2}\int_\T\partial_i|\nabla V|^2\dd x=0.
\]
The other terms in \eqref{eq:momentum} are divergences, and
$\int_\T(m-\rho u)\dd x=\int_\T\nabla\rho\dd x=0$, so
$\int_\T\rho u\dd x$ and $\int_\T\rho v\dd x$ are conserved and equal. With
$\int_\T\rho\dd x=M$, both equal $MU$; centering either kinetic energy at $U$
therefore subtracts the constant $M|U|^2/2$.

For $\psi=u,v$, the continuity equation gives
\[
\int_\Omega\rho\psi\cdot(\psi_t+u\cdot\nabla\psi)\dd x
=\frac d{dt}\int_\Omega\tfrac12\rho|\psi|^2\dd x
-\frac12\int_\Omega(\rho_t+\Div(\rho u))|\psi|^2\dd x
=\frac d{dt}\int_\Omega\tfrac12\rho|\psi|^2\dd x.
\]
On each compact interval the local construction in Appendix~\ref{app:local}
gives $u,v\in C H^2$, $u_t,v_t\in C L^2$, $s-\sqrt b\in C H^3$ and
$s_t\in C H^1$. In $\R^3$ the transport flux $\rho u|\psi|^2$ is integrable
because $\rho u\in L^2$ and $\psi\in L^4$; spatial cutoffs therefore justify
the integration by parts at infinity.

For the internal and electric energies, the chain rule and $QE=E$ give
\begin{align}
\frac d{dt}\int_\Omega H(\rho)\dd x
&=\int_\Omega H'(\rho)\rho_t\dd x
=-\int_\Omega H'(\rho)\Div(\rho u)\dd x
=\int_\Omega\rho u\cdot\nabla H'(\rho)\dd x
=\int_\Omega u\cdot\nabla\Pi\dd x,\label{eq:internal-chain}\\
\frac d{dt}\frac{\ell^2}{2}\|E\|_2^2
&=\ell^2\int_\Omega E\cdot E_t\dd x
=\int_\Omega E\cdot Q(\rho u)\dd x
=\int_\Omega\rho u\cdot E\dd x,\label{eq:field-chain}
\end{align}
where we used $H''(\rho)=P'(\rho)/\rho$, so that
$\rho\nabla H'(\rho)=P'(\rho)\nabla\rho=\nabla\Pi$, and
$E_t=\ell^{-2}Q(\rho u)$. For the quantum term,
$\Div(\rho\nabla^2\log\rho)=2\rho\nabla(\Delta s/s)$ and \eqref{eq:mass} give
\begin{equation}\label{eq:quantum-work}
\begin{aligned}
\int_\Omega u\cdot\Div(\rho\nabla^2\log\rho)\dd x
&=2\int_\Omega\rho u\cdot\nabla\Bigl(\frac{\Delta s}{s}\Bigr)\dd x
=-2\int_\Omega\Div(\rho u)\frac{\Delta s}{s}\dd x\\
&=2\int_\Omega\rho_t\frac{\Delta s}{s}\dd x
=4\int_\Omega s_t\Delta s\dd x
=-2\frac d{dt}\|\nabla s\|_2^2.
\end{aligned}
\end{equation}
Testing \eqref{eq:momentum} with $u$ and using
\[
\int_\Omega u\cdot\bigl((\rho u)_t+\Div(\rho u\otimes u)\bigr)\dd x
=\frac d{dt}\int_\Omega\tfrac12\rho|u|^2\dd x
\]
together with \eqref{eq:internal-chain}--\eqref{eq:quantum-work} gives
\eqref{eq:Eu}.

Similarly, testing \eqref{eq:effective-velocity} with $v$ and using
\eqref{eq:internal-chain}--\eqref{eq:field-chain}, the terms containing $u$
reproduce the derivatives of the internal and field energies, so
\[
E_v'+\int_\Omega\rho|\nabla v|^2\dd x
=-\int_\Omega\nabla\log\rho\cdot\nabla\Pi\dd x
-\int_\Omega\nabla r\cdot E\dd x.
\]
By $\nabla\Pi=\gamma\rho^{\gamma-1}\nabla\rho$, the first term on the right is
$-\int_\Omega\nabla\log\rho\cdot\nabla\Pi\dd x=-\gamma\int_\Omega\rho^{\gamma-2}|\nabla\rho|^2\dd x$.
Since $\Div E=-\ell^{-2}r$, the second term is
\[
-\int_\Omega\nabla r\cdot E\dd x=\int_\Omega r\,\Div E\dd x=-\ell^{-2}\|r\|_2^2,
\]
which gives \eqref{eq:Ev}.
\end{proof}

To estimate concentration we use a Hessian bound whose constant is independent
of the minimum density. In the lemma, $\mathcal I$ is evaluated at a fixed time.
\begin{lemma}\label{lem:quantum}
Let $r\in H^3(\Omega)$ and $\inf_\Omega\rho>0$. Write
$\mathcal I=\int_\Omega\rho|\nabla^2\log\rho|^2\dd x$. Then
\begin{equation}\label{eq:quantum}
   \|\nabla^2 s\|_2\le (1+\sqrt3/4) \mathcal I^{1/2}.
\end{equation}
On a compact interval $[0,T]$ of positive strong existence,
\begin{equation}\label{eq:Iintegral}
   \int_0^T \mathcal I(t)\dd t\le E_u(0)+2E_v(0).
\end{equation}
On the unit torus one also has $\|s(t)\|_{H^1}^2\le M+E_u(0)/2$.
\end{lemma}
\begin{proof}
With $g=\log\rho$, the divergence identity
\[
\Div(\rho|\nabla g|^2\nabla g)
=\rho|\nabla g|^4+2\rho\,\nabla^2g:(\nabla g\otimes\nabla g)
+\rho|\nabla g|^2\Delta g
\]
holds almost everywhere. On the torus the divergence integrates to zero by periodicity.
In $\R^3$, $r\in H^3$ and $\inf\rho>0$ give $\rho\in L^\infty$,
$\nabla g\in L^2\cap L^\infty$ and $\nabla^2g\in L^2$. Thus all three terms
on the right belong to $L^1$, and the flux $F=\rho|\nabla g|^2\nabla g$ satisfies
\[
 \|F\|_1\le\|\rho\|_\infty\|\nabla g\|_\infty\|\nabla g\|_2^2.
\]
Choose $0\le\chi\le1$ in $C_c^\infty(\R^3)$ with $\chi=1$ on the unit ball,
and set $\chi_N(x)=\chi(x/N)$. Testing the identity with $\chi_N$ gives
$\bigl|\int F\cdot\nabla\chi_N\dd x\bigr|\le N^{-1}\|\nabla\chi\|_\infty\|F\|_1$.
Letting $N\to\infty$ gives zero for the integral of the divergence. Since
$|\nabla^2g:(\nabla g\otimes\nabla g)|\le|\nabla^2g||\nabla g|^2$ and
$|\Delta g|\le\sqrt3|\nabla^2g|$, Cauchy--Schwarz gives
\[
\int_\Omega\rho|\nabla g|^4\dd x
\le(2+\sqrt3)\int_\Omega\rho|\nabla^2g||\nabla g|^2\dd x
\le(2+\sqrt3)\,\mathcal I^{1/2}
 \Bigl(\int_\Omega\rho|\nabla g|^4\dd x\Bigr)^{1/2}.
\]
Now $\nabla^2s=\tfrac12s\nabla^2g+\tfrac14s\nabla g\otimes\nabla g$, hence
\[
\|\nabla^2s\|_2
\le\tfrac12\mathcal I^{1/2}
+\tfrac14\Bigl(\int_\Omega\rho|\nabla g|^4\dd x\Bigr)^{1/2}
\le\Bigl(\tfrac12+\tfrac{2+\sqrt3}{4}\Bigr)\mathcal I^{1/2}
=\Bigl(1+\tfrac{\sqrt3}{4}\Bigr)\mathcal I^{1/2},
\]
which is \eqref{eq:quantum}. Finally, $\nabla^2g=Dv-Du$, so
$|\nabla^2g|^2=|Dv-Du|^2\le2(|\nabla v|^2+|Du|^2)$ and
$\mathcal I\le2\int_\Omega\rho(|\nabla v|^2+|Du|^2)\dd x$. Integration of
\eqref{eq:Eu}--\eqref{eq:Ev} proves \eqref{eq:Iintegral}. The periodic $H^1$
bound follows from $\|s\|_2^2=M$, the inequality $2\|\nabla s\|_2^2\le E_u$ and
the monotonicity of $E_u$.
\end{proof}

Let $T^*$ be the lifespan of the maximal positive strong solution
given by Proposition~\ref{prop:local-existence}. On each compact interval
$[0,S]\subset[0,T^*)$, the continuity and far-field argument following
Definition~\ref{def:positive-strong-solution} gives
\begin{equation}\label{eq:compact-positive}
  0<\inf_{[0,S]\times\Omega}\rho\le\sup_{[0,S]\times\Omega}\rho<\infty.
\end{equation}
The estimates below quantify these bounds independently of the minimum
on the interval of existence. Write $R(t)=\|\rho(t)\|_\infty$ and
$Y_p(t)=\int_\Omega\rho|v|^p\dd x$ for $p\ge2$. We fix constants
$S_1,S_2\ge1$ for the inequalities $\|g\|_6\le S_1\|\nabla g\|_2$ on
$\R^3$ and $\|g\|_\infty\le S_2\|g\|_{H^2}$ on both domains, also for
vector- and tensor-valued functions. In the first inequality $\dot H^1$
is realized by its $L^6$ representatives.

\subsection{Periodic elliptic and interpolation estimates}
On the unit torus we denote spatial means by an overbar and take $b=M$.
The elliptic equation \eqref{eq:poisson} will be used both in Sobolev spaces and in pointwise estimates. The required bounds follow from its Fourier multiplier and the periodic heat kernel.
\begin{lemma}\label{lem:poisson-estimates}
If $\int f\dd x=0$ and $-\ell^2\Delta V=f$, $\int V\dd x=0$, then
\[
 \|V\|_{H^{s+2}}\le C_s\ell^{-2}\|f\|_{H^s}\quad(s\ge0).
\]
If $1\le p_1\le p_2\le\infty$ and $1/p_1-1/p_2<1/3$, then
\begin{equation}\label{eq:poisson-Lp}
 \|\nabla V\|_{p_2}\le C_{p_1,p_2}\ell^{-2}\|f\|_{p_1}.
\end{equation}
\end{lemma}
\begin{proof}
At every nonzero Fourier frequency, $\widehat V(\xi)=\ell^{-2}|\xi|^{-2}\widehat f(\xi)$. The Sobolev bound follows from Parseval's identity and $|\xi|\ge2\pi$. To estimate the gradient, write
\[
 \nabla V=\ell^{-2}\int_0^\infty\nabla e^{t\Delta}f\dd t.
\]
The periodized Gaussian derivative and Young's inequality bound the
integrand in $L^{p_2}$ by
$Ct^{-1/2-(3/2)(1/p_1-1/p_2)}\|f\|_{p_1}$ for $0<t\le1$.
For $t\ge1$, the differentiated periodic heat kernel has exponentially
decaying Fourier series, giving $Ce^{-ct}\|f\|_{p_1}$ instead.
Since $1/2+(3/2)(1/p_1-1/p_2)<1$, integration in time gives \eqref{eq:poisson-Lp}.
This also justifies the representation, first for smooth $f$
and then by approximation (or weak limits when $p_1=\infty$).
\end{proof}

The next interpolation estimate turns \eqref{eq:quantum} and \eqref{eq:Iintegral} into the time integrability needed at cubic pressure.

\begin{lemma}\label{lem:agmon}
For $f\in H^2(\T)$,
\begin{equation}\label{eq:agmon}
 \|f\|_\infty^2\le C\|f\|_{H^1}\|f\|_{H^2}.
\end{equation}
Consequently, for each finite $L$ and every $T<\min\{T^*,L\}$,
\begin{equation}\label{eq:rhoL2Linf}
 \sup_{0\le t\le T}\|\rho(t)\|_3\le C,
 \qquad \int_0^T\|\rho(t)\|_\infty^2\dd t\le C(1+L).
\end{equation}
The constant $C$ depends only on $M$, the two initial energies, and the unit torus.
\end{lemma}
\begin{proof}
The lattice bounds $\sum_{0<|\xi|\le R}|\xi|^{-2}\le CR$ and
$\sum_{|\xi|>R}|\xi|^{-4}\le C/R$ hold for $R\ge1$. Splitting the Fourier series at $R$ and applying Cauchy--Schwarz gives
\[
 \|f\|_\infty\le C\bigl(R^{1/2}\|f\|_{H^1}+R^{-1/2}\|f\|_{H^2}\bigr).
\]
For $f\ne0$, choose $R=\|f\|_{H^2}/\|f\|_{H^1}$ to obtain \eqref{eq:agmon}. For $f=\sqrt\rho$, \eqref{eq:quantum}--\eqref{eq:Iintegral} thus imply
\[
 \|\rho(t)\|_\infty^2\le C\|\sqrt\rho\|_{H^1}^2\|\sqrt\rho\|_{H^2}^2
 \le C(1+\mathcal I(t)),\qquad \|\rho(t)\|_3=\|\sqrt\rho\|_6^2\le C.
\]
Integration proves \eqref{eq:rhoL2Linf}.
\end{proof}

\subsection{Dissipation and field estimates in the whole space}
In all whole-space estimates below, write $\budget=1+E_u(0)+E_v(0)$, with $U=0$.
In this subsection take $\Omega=\R^3$ and set $f=\sqrt\rho-\sqrt b$. Define
\[
 h_{b,\gamma}=\inf_{z>0,\ z\ne b}\frac{H(z\mid b)}{(\sqrt z-\sqrt b)^2}.
\]
The following estimates account for the positive background in the density
energy and retain the charge dissipation needed for uniform bounds at late times.
\begin{proposition}\label{prop:energy-bounds}
Suppose that the data satisfy the hypotheses of Theorem~\ref{thm:cauchy}
and that a positive strong solution exists on $[0,S]$. Then
$h_{b,\gamma}=b^{\gamma-1}$ and
\begin{equation}\label{eq:charge-budget}
 \sup_{0\le t\le S}\|r(t)\|_2\le(\budget/h_{b,\gamma})^{1/4}
 (S_1\sqrt{\budget/2})^{3/2}+2\sqrt b\,(\budget/h_{b,\gamma})^{1/2}.
\end{equation}
Moreover,
\begin{equation}\label{eq:dissipation-budget}
 \int_0^S\Bigl(\|r\|_2^2+\|\nabla^2s\|_2^2+
 \int\rho(|\nabla v|^2+|\nabla u|^2)\dd x\Bigr)\dd t
 \le(7+\ell^2)\budget.
\end{equation}
Every subinterval $I\subset[0,S]$ satisfies
\begin{equation}\label{eq:density-time-budget}
 \int_I R(t)^2\dd t\le8b^2|I|+\frac{80}{\pi^4}\budget^2.
\end{equation}
\end{proposition}
\begin{proof}
The initial far-field condition implies $c_0\le b\le C_0$. Since
$v_0=u_0+\nabla r_0/\rho_0$, we have
$|v_0|^2\le2|u_0|^2+2|\nabla r_0|^2/\rho_0^2$, and since
$2|\nabla\sqrt{\rho_0}|^2=|\nabla r_0|^2/(2\rho_0)$,
\[
\tfrac12\rho_0(|u_0|^2+|v_0|^2)+2|\nabla\sqrt{\rho_0}|^2
\le\tfrac32\rho_0|u_0|^2+\tfrac3{2\rho_0}|\nabla r_0|^2.
\]
Taylor's formula gives
$2H(\rho_0)\le r_0^2\max_{c_0\le z\le C_0}\gamma z^{\gamma-2}$.
Thus
\begin{align}
 \budget\le 1&+\tfrac32 C_0\|u_0\|_2^2+\tfrac3{2c_0}\|\nabla r_0\|_2^2
 +\ell^2\|E_0\|_2^2\notag\\
 &+\Bigl(\max_{c_0\le z\le C_0}\gamma z^{\gamma-2}\Bigr)\|r_0\|_2^2.
 \label{eq:initial-energy}
\end{align}

To identify $h_{b,\gamma}$, put $y=\sqrt{z/b}$. For $\gamma>1$,
\[
 H(z)-b^{\gamma-1}(\sqrt z-\sqrt b)^2
 =\frac{b^\gamma y}{\gamma-1}
    \bigl[y^{2\gamma-1}-1-(2\gamma-1)(y-1)\bigr]\ge0,
\]
since the tangent at $y=1$ of the convex function $y\mapsto y^{2\gamma-1}$
is $y\mapsto1+(2\gamma-1)(y-1)$. For $\gamma=1$ the same difference is
$2by(y\log y-y+1)\ge0$. The quotient defining $h_{b,\gamma}$ tends to
$b^{\gamma-1}$ as $z\downarrow0$.

Since $E_u$ and $E_v$ are nonincreasing by \eqref{eq:Eu}--\eqref{eq:Ev}, the
identities $H(\rho)\ge h_{b,\gamma}f^2$ and \eqref{eq:energy-u}--\eqref{eq:energy-v}
give
\begin{equation}\label{eq:basic-bounds}
 Y_2\le2\budget,\quad \|E\|_2^2\le2\ell^{-2}\budget,\quad
 \|f\|_2^2\le\budget/h_{b,\gamma},\quad
 \|\nabla f\|_2^2\le\budget/2.
\end{equation}
Indeed, $Y_2\le2E_v(0)$, $\ell^2\|E\|_2^2/2\le E_v(0)$,
$h_{b,\gamma}\|f\|_2^2\le E_v(0)$ and $2\|\nabla f\|_2^2\le E_u(0)$, and
$\budget=1+E_u(0)+E_v(0)$.
Since $r=f^2+2\sqrt b\,f$,
\begin{align*}
 \|r\|_2&\le\|f\|_4^2+2\sqrt b\,\|f\|_2,\\
 \|f\|_4^2&\le\|f\|_2^{1/2}\|f\|_6^{3/2}
 \le\|f\|_2^{1/2}(S_1\|\nabla f\|_2)^{3/2},
\end{align*}
which gives \eqref{eq:charge-budget}. Interpolation for the field,
$\|E\|_4^4\le\|E\|_2\|E\|_6^3$, and \eqref{eq:gauss-symbol} also give
\begin{equation}\label{eq:field-four}
 \|E\|_4^4
 \le\ell^{-1}\sqrt{2\budget}\,(S_1\ell^{-2}\|r\|_2)^3.
\end{equation}

We shall also use a pointwise estimate for $E$. For smooth compactly
supported $r$, the field is the convolution of $r$ with
$-x/(4\pi\ell^2|x|^3)$. Estimating the integral separately on the unit
ball and its complement gives
\begin{align}
 \|E\|_\infty
 &\le\frac1{4\pi\ell^2}
 \left(\|r\|_\infty\int_{|y|\le1}|y|^{-2}\dd y
 +\|r\|_2\Bigl(\int_{|y|>1}|y|^{-4}\dd y\Bigr)^{1/2}\right)\notag\\
 &=\ell^{-2}\bigl(\|r\|_\infty+(4\pi)^{-1/2}\|r\|_2\bigr).
 \label{eq:newton-field}
\end{align}
For general $r$ at a fixed time, approximate in $H^3$ by
$r_n\in C_c^\infty$. The corresponding Newton fields $E_n$ belong to
$L^2$ since they are bounded and decay as $O(|x|^{-2})$.
By \eqref{eq:gauss-symbol},
$\|E_n-E\|_6\le S_1\ell^{-2}\|r_n-r\|_2\to0$. Along an almost
everywhere convergent subsequence we may pass to the bound
\eqref{eq:newton-field}, since also $r_n\to r$ in $L^\infty$.

The Hessian bound \eqref{eq:quantum} applies also to $\Omega=\R^3$.
Writing $g=\log\rho$, we have $\nabla v-\nabla u=\nabla^2g$, a symmetric matrix. Thus
\[
 |\nabla^2g|^2=|Dv-Du|^2\le2(|\nabla v|^2+|Du|^2),\qquad
 |\nabla u|^2=|Du|^2+|\nabla v-Dv|^2\le |Du|^2+|\nabla v|^2.
\]
Since $2(1+\sqrt3/4)^2<5$, the energy identities yield
\begin{equation}\label{eq:hessian-time}
 \int_0^S\|\nabla^2s\|_2^2\dd t
 \le 5\int_0^S\!\int\rho(|\nabla v|^2+|Du|^2)\dd x\dd t
 \le5\bigl(E_v(0)+\tfrac12E_u(0)\bigr)\le5\budget.
\end{equation}
Adding the two weighted velocity gradients and the charge dissipation
bounds the left-hand side of \eqref{eq:dissipation-budget} by
$(\ell^2+7)E_v(0)+3E_u(0)\le(\ell^2+7)\budget$; here we used
$|\nabla u|^2\le|Du|^2+|\nabla v|^2$, so that
$\int\!\!\int\rho|\nabla u|^2\le\frac12E_u(0)+E_v(0)$.

Finally, splitting Fourier inversion at radius $\lambda>0$ and applying
Cauchy--Schwarz on each side gives
\begin{equation}\label{eq:frequency-split}
 (2\pi)^{3/2}\|f\|_\infty
 \le (4\pi\lambda)^{1/2}\|\nabla f\|_2+
 (4\pi/\lambda)^{1/2}\|\nabla^2s\|_2 ,
\end{equation}
because
\[
\int_{|\xi|\le\lambda}|\xi|^{-2}\dd\xi=4\pi\lambda,
\qquad
\int_{|\xi|>\lambda}|\xi|^{-4}\dd\xi=\frac{4\pi}{\lambda}.
\]
Choosing $\lambda=\|\nabla^2s\|_2/\|\nabla f\|_2$ minimizes the right-hand
side, the minimum being
\[
2(4\pi)^{1/2}\|\nabla f\|_2^{1/2}\|\nabla^2s\|_2^{1/2}.
\]
Hence
\begin{equation}\label{eq:fourier-density}
 \|f\|_\infty^2\le\frac{2}{\pi^2}\|\nabla f\|_2\|\nabla^2s\|_2.
\end{equation}
If either derivative norm vanishes, then $f=0$ since $f\in L^2$.
Now $R^2\le(\sqrt b+\|f\|_\infty)^4\le8b^2+8\|f\|_\infty^4$, and
\eqref{eq:fourier-density} gives
\[
R^2\le8b^2+\frac{32}{\pi^4}\|\nabla f\|_2^2\|\nabla^2s\|_2^2.
\]
Therefore \eqref{eq:basic-bounds} and \eqref{eq:hessian-time} yield
\[
 \int_I R^2\dd t\le8b^2|I|+
 \frac{32}{\pi^4}\sup_I\|\nabla f\|_2^2
             \int_I\|\nabla^2s\|_2^2\dd t
 \le8b^2|I|+\frac{80}{\pi^4}\budget^2.
\]
\end{proof}

\section{Upper bounds for the density}\label{sec:upper-both}
The time integral of $R^2$ controls the pressure work up to $\gamma=3$.
The periodic calculation uses the finite mass to bound a mixed moment.
In the whole space we subtract the far-field contribution before taking
the density moment.

\subsection{Heat smoothing and mixed moments}
Fix a finite horizon $L>0$. In this subsection the periodic solution is considered
on $[0,T]$, where $T<\min\{T^*,L\}$. All constants are independent of $T$
and of the minimum of the density on this interval. The two inputs from
Lemma~\ref{lem:agmon} are
\begin{equation}\label{eq:inputs}
 \sup_{0\le t\le T}\|\rho(t)\|_3\le C,
 \qquad \int_0^T\|\rho(t)\|_\infty^2\dd t\le C_L.
\end{equation}
The constants in \eqref{eq:inputs} depend only on the initial energies
and mass, with additional dependence on $L$ in the second bound.

Through the continuity equation \eqref{eq:heat-mass}, a bound for the effective
momentum $m=\rho v$ in $L^\infty(0,T;L^p)$ with $p>3$ gives a pointwise
upper bound for the density on $[0,T]$.

\begin{lemma}\label{lem:heat}
Let $\Omega=\T$ and $b=M$. Suppose that $p>3$, $0<T\le L$ and
$\sup_{0\le t\le T}\|\rho v(t)\|_p\le B_p$. Then
\begin{equation}\label{eq:heat-bound}
 \sup_{0\le t\le T}\|\rho(t)\|_\infty
 \le \|\rho_0\|_\infty+
 C_pB_p\bigl(L+L^{\frac12-\frac{3}{2p}}\bigr).
\end{equation}
\end{lemma}
\begin{proof}
The periodic heat kernel $\mathcal K_h$ is the periodization of the Euclidean
Gaussian at time $h>0$. Its first derivative satisfies
\begin{equation}\label{eq:periodic-kernel}
\|\nabla \mathcal K_h\|_{p'}\le C_ph^{-1/2-3/(2p)}
\qquad(0<h\le1,\ 1/p+1/p'=1).
\end{equation}
Indeed, this follows by scaling the Gaussian derivative on $\R^3$ and summing
its exponentially decaying translates. For $h\ge1$, the semigroup identity
and the $L^1$ contraction give a fixed bound for the same norm. Young's
convolution inequality consequently gives
\[
 \|\Div e^{h\Delta}\mathbf f\|_\infty
 \le C_p\bigl(1+h^{-1/2-3/(2p)}\bigr)\|\mathbf f\|_p.
\]
Using \eqref{eq:heat-mass}, write
\[
 \rho(t)=e^{t\Delta}\rho_0-
 \int_0^t\Div e^{(t-\tau)\Delta}(\rho v)(\tau)\dd\tau.
\]
Since $\alpha_p=1/2+3/(2p)<1$, contraction in $L^\infty$ and the integration
$\int_0^t(1+(t-\tau)^{-\alpha_p})\dd\tau
=t+t^{1-\alpha_p}/(1-\alpha_p)$ with $1-\alpha_p=1/2-3/(2p)$ give
\eqref{eq:heat-bound}.
\end{proof}

\paragraph{Polytropic pressure.}
Let $1<\gamma\le3$ and $a>1$. Consider the conservative state
$W=(\rho,m)\in(0,\infty)\times\R^3$. In state differentiation,
$\rho$ and $m$ are independent coordinates. Along a solution, the
functions below are evaluated at $W=W(x,t)$. Define
\begin{equation}\label{eq:tm-entropy}
 \begin{gathered}
 v=m/\rho,\qquad
 \mathfrak e(W)=1+\frac{|m|^2}{2\rho^2}+\frac{\rho^{\gamma-1}}{\gamma-1},\\
 \Psi_a(W)=\rho \mathfrak e(W)^a.
 \end{gathered}
\end{equation}
We abbreviate $\mathfrak e(W)$ to $\mathfrak e$ in the calculations below. The lower bound $\mathfrak e\ge1$ will be used in estimating the electric work.

The convexity of $\Psi_a$ controls diffusion in the conservative variables. The following identity keeps the nonzero transport contribution explicit.

\begin{lemma}\label{lem:mixed}
In the periodic setting with $b=M$, let $1<\gamma\le3$ and $a>1$. Put
$Z_a=\int\rho \mathfrak e^a\dd x$, with $\mathfrak e$ evaluated
at $(\rho,m)$. Then
\begin{align}
 Z_a'&+a\int\rho \mathfrak e^{a-1}|\nabla v|^2\dd x
       +a\gamma\int\rho^{\gamma-2}\mathfrak e^{a-1}|\nabla\rho|^2\dd x\notag\\
 &\quad+a(a-1)\int\rho \mathfrak e^{a-2}|\nabla \mathfrak e|^2\dd x\notag\\
 &=a(a-1)\int\rho^\gamma \mathfrak e^{a-2}v\cdot\nabla \mathfrak e\dd x
   -a\int\rho \mathfrak e^{a-1}v\cdot\nabla V\dd x.\label{eq:tm-mixed-id}
\end{align}
In particular,
\begin{equation}\label{eq:tm-growth}
 Z_a'\le C_{a,\gamma,\ell,M}
 \bigl(1+\|\rho\|_\infty^{\gamma-1}+\|\rho\|_\infty\bigr)Z_a.
\end{equation}
\end{lemma}
\begin{proof}
Fix a direction $\zeta=(\zeta_\rho,\zeta_m)\in\R\times\R^3$
and differentiate along the affine line $W+s\zeta$ at $s=0$.
The induced first variations of $v$ and $\mathfrak e$ are
\[
 \zeta_v:=D_Wv(W)[\zeta]=\frac{\zeta_m-v\zeta_\rho}{\rho},
 \qquad
 \zeta_{\mathfrak e}:=D_W\mathfrak e(W)[\zeta]
       =v\cdot\zeta_v+\rho^{\gamma-2}\zeta_\rho.
\]
The function $\rho \mathfrak e=\rho+|m|^2/(2\rho)+\rho^\gamma/(\gamma-1)$ has second variation $\rho|\zeta_v|^2+\gamma\rho^{\gamma-2}\zeta_\rho^2$. As the first state coordinate is affine along the line,
\begin{align}
 D_W^2\Psi_a(W)[\zeta,\zeta]
 &=2a\mathfrak e^{a-1}\zeta_\rho\zeta_{\mathfrak e}
     +a\rho \mathfrak e^{a-1}D_W^2\mathfrak e(W)[\zeta,\zeta]
     +a(a-1)\rho \mathfrak e^{a-2}\zeta_{\mathfrak e}^2\nonumber\\
 &=a \mathfrak e^{a-1}D_W^2(\rho \mathfrak e)(W)[\zeta,\zeta]
    +a(a-1)\rho \mathfrak e^{a-2}\zeta_{\mathfrak e}^2\nonumber\\
 &=a\rho \mathfrak e^{a-1}|\zeta_v|^2
   +a\gamma\rho^{\gamma-2}\mathfrak e^{a-1}\zeta_\rho^2
   +a(a-1)\rho \mathfrak e^{a-2}\zeta_{\mathfrak e}^2.\label{eq:tm-hessian}
\end{align}
In particular $\Psi_a$ is convex. Testing the diffusion part $\Delta W$
against $\nabla_W\Psi_a$ and integrating by parts gives
$-\sum_j\int D_W^2\Psi_a[\partial_jW,\partial_jW]\dd x$, which supplies the
three dissipations in \eqref{eq:tm-mixed-id}.
To compute the transport contribution, drop the Laplacian in
\eqref{eq:heat-mass}--\eqref{eq:heat-momentum} and the electric term. The
resulting first-order subsystem is transported by $v$, and its equations are
\begin{align*}
 (\partial_t+v\cdot\nabla)\rho&=-\rho\Div v,\\
 \rho(\partial_t+v\cdot\nabla)v&=-\nabla\rho^\gamma,\\
 \rho(\partial_t+v\cdot\nabla)\mathfrak e
 &=-v\cdot\nabla\rho^\gamma-\rho^\gamma\Div v
  =-\Div(\rho^\gamma v).
\end{align*}
Consequently,
\[
 \partial_t(\rho \mathfrak e^a)+\Div(\rho v\mathfrak e^a)
 =-a\mathfrak e^{a-1}\Div(\rho^\gamma v).
\]
The two remaining contributions to the integrated identity are therefore
\begin{align*}
 -a\int \mathfrak e^{a-1}\Div(\rho^\gamma v)\dd x
 &=a\int\rho^\gamma v\cdot\nabla(\mathfrak e^{a-1})\dd x\\
 &=a(a-1)\int\rho^\gamma \mathfrak e^{a-2}v\cdot\nabla \mathfrak e\dd x,\\
 \nabla_m\Psi_a(W)&=a\rho \mathfrak e^{a-1}\frac m{\rho^2}
                    =a\mathfrak e^{a-1}v,\\
 \int\nabla_m\Psi_a(W)\cdot(-\rho\nabla V)\dd x
 &=-a\int\rho \mathfrak e^{a-1}v\cdot\nabla V\dd x.
\end{align*}
Here $\nabla_m$ differentiates at fixed $\rho$. Adding the diffusion contribution gives \eqref{eq:tm-mixed-id}. By Young's inequality the transport remainder is bounded by
\[
 \tfrac12 a(a-1)\int \rho \mathfrak e^{a-2}|\nabla \mathfrak e|^2\dd x
 + C_a\int\rho^{2\gamma-1}\mathfrak e^{a-2}|v|^2\dd x.
\]
Since $|v|^2\le2\mathfrak e$ and $\rho^{\gamma-1}\le(\gamma-1)\mathfrak e$, the second integral is at most $2(\gamma-1)\|\rho\|_\infty^{\gamma-1} Z_a$. Also $\mathfrak e\ge1$ gives
\[
 \left|a\int\rho \mathfrak e^{a-1}v\cdot\nabla V\dd x\right|
 \le \sqrt2a\|\nabla V\|_\infty Z_a,\qquad
 \|\nabla V\|_\infty\le C_\ell(\|\rho\|_\infty+M),
\]
the last bound following from Lemma~\ref{lem:poisson-estimates} with
$p_1=4$, the interpolation $\|\rho\|_4\le M^{1/4}\|\rho\|_\infty^{3/4}$ and
Young's inequality
$M^{1/4}\|\rho\|_\infty^{3/4}\le\tfrac14M+\tfrac34\|\rho\|_\infty$.
Absorbing half the $|\nabla \mathfrak e|^2$ term proves \eqref{eq:tm-growth}.
\end{proof}

The coefficient in \eqref{eq:tm-growth} is integrable by \eqref{eq:inputs}. This gives a fixed mixed moment large enough to apply Lemma~\ref{lem:heat}.

\begin{proposition}\label{prop:upper-poly}
Let $\Omega=\T$ and $b=M$. If $1<\gamma\le3$, there is a constant $R_L$,
depending only on $L$, the fixed parameters, the initial quantities in
Theorem~\ref{thm:periodic} and absolute constants, such that
\begin{equation}\label{eq:tm-upper}
 \sup_{0\le t\le T}\|\rho(t)\|_\infty\le R_L
 \qquad\bigl(T<\min\{T^*,L\}\bigr).
\end{equation}
\end{proposition}
\begin{proof}
Since $0<\gamma-1\le2$, we have $R^{\gamma-1}\le1+R^2$, where
$R=\|\rho\|_\infty$. With $C_L$ as in \eqref{eq:inputs}, Cauchy--Schwarz gives
$\int_0^T R\dd t\le\sqrt{LC_L}$. The time integral of the coefficient in
\eqref{eq:tm-growth} is therefore bounded by
\[
 K_{a,L}:=C_{a,\gamma,\ell,M}\bigl(2L+C_L+\sqrt{LC_L}\bigr).
\]
Gronwall's inequality gives $\sup_{t\le T}Z_a(t)\le Z_a(0)e^{K_{a,L}}$.
The initial value is controlled by
$\rho_0,v_0\in L^\infty$ and the mass, since
$v_0=u_0+\nabla\log\rho_0\in H^2\hookrightarrow L^\infty$
under \eqref{eq:initial-data}. More explicitly,
\begin{align*}
 Z_a(0)&\le M\left(1+\tfrac12\|v_0\|_\infty^2
                 +\frac{\|\rho_0\|_\infty^{\gamma-1}}{\gamma-1}\right)^a,\\
 \int\rho|v|^{2a}\dd x&\le2^a\int\rho \mathfrak e^a\dd x=2^aZ_a(t)\le C_{a,L},\\
 \int\rho^{1+a(\gamma-1)}\dd x
 &\le(\gamma-1)^a\int\rho \mathfrak e^a\dd x
   =(\gamma-1)^aZ_a(t)\le C_{a,L},
\end{align*}
where we used $|v|^{2a}\le2^a\mathfrak e^a$ and
$\rho^{\gamma-1}\le(\gamma-1)\mathfrak e$. Set $p=2a$ and
$q_d=1+a(\gamma-1)$. Factoring $\rho v=(\rho^{1/p}v)\rho^{1-1/p}$ and
applying H\"older's inequality gives
\begin{equation}\label{eq:tm-flux}
 \|\rho v\|_{p_*}
 =\|(\rho^{1/p}v)\rho^{1-1/p}\|_{p_*}
 \le\Bigl(\int\rho|v|^p\dd x\Bigr)^{1/p}\|\rho\|_{q_d}^{1-1/p},
\end{equation}
with
\[
 \frac1{p_*}=\frac1p+\frac{1-1/p}{q_d}
 =\frac{1}{2a}+\frac{1-1/(2a)}{1+a(\gamma-1)}
 =\frac{\gamma+1}{2[1+a(\gamma-1)]},
\]
that is, $p_*=2[1+a(\gamma-1)]/(\gamma+1)$. Choose once and for all
$a>(3\gamma+1)/(2(\gamma-1))$. Then
\begin{align*}
 p_*-3&=\frac{2a(\gamma-1)-3\gamma-1}{\gamma+1}>0,\\
 \sup_{0\le t\le T}\|\rho(t)\|_\infty
 &\le\|\rho_0\|_\infty+C_{p_*}C_{a,L}
            \left(L+L^{\frac12-\frac3{2p_*}}\right)=:R_L,
\end{align*}
where the second line is the integral estimate of Lemma~\ref{lem:heat}
with the flux just obtained. At $\gamma=3$ one may take
$a=3$, so that $(p,q_d,p_*)=(6,7,7/2)$ and the heat-kernel singularity
is $1/2+3/(2p_*)=13/14<1$. \end{proof}

\paragraph{Isothermal pressure.}
At $\gamma=1$ the internal-energy weight in \eqref{eq:tm-entropy} is singular. A twelfth velocity moment and a fourth density moment yield the required flux integrability.

\begin{proposition}\label{prop:upper-iso}
Under the hypotheses of Theorem~\ref{thm:periodic} with $\gamma=1$,
the estimate \eqref{eq:tm-upper} holds with the same permitted dependence
of its constant.
\end{proposition}
\begin{proof}
First, \eqref{eq:inputs} and Lemma~\ref{lem:poisson-estimates} with
$p_1=3$, $p_2=18$ give $\|\nabla V\|_{18}\le C_\ell\|\rho-M\|_3\le C$.
Here $1/3-1/18=5/18<1/3$, so the estimate is strictly away from its
non-integrable endpoint. Hence
\begin{equation}\label{eq:tm-field12}
 \begin{aligned}
 \left(\int\rho|\nabla V|^{12}\dd x\right)^{1/12}
 &\le\left[\left(\int\rho^3\dd x\right)^{1/3}
             \left(\int|\nabla V|^{18}\dd x\right)^{2/3}\right]^{1/12}\\
 &=\|\rho\|_3^{1/12}\|\nabla V\|_{18}
 \le C_\ell\|\rho\|_3^{1/12}(\|\rho\|_3+M)\le C.
 \end{aligned}
\end{equation}
Let $Y_{12}=\int\rho|v|^{12}\dd x$. Test
\eqref{eq:effective-velocity} by $|v|^{10}v$, using \eqref{eq:mass}
to cancel the transport with velocity $u$. Since $\rho^\gamma=\rho$,
\begin{align}
 \frac1{12}Y_{12}'
 &+\int\rho|v|^{10}|\nabla v|^2\dd x
  +10\int\rho|v|^{10}|\nabla|v||^2\dd x\nonumber\\
 &=\int\rho\Div(|v|^{10}v)\dd x
   -\int\rho|v|^{10}v\cdot\nabla V\dd x.\label{eq:tm-moment12}
\end{align}
Here the last diffusion integrand may equivalently be written
$\rho|v|^8\sum_j(v\cdot\partial_jv)^2$, with value zero when $v=0$.
The pressure term is bounded using $|\Div(|v|^{10}v)|\le(\sqrt3+10)|v|^{10}|\nabla v|$. Young's inequality and conservation of mass give
\[
 \left|\int\rho\Div(|v|^{10}v)\dd x\right|
 \le \tfrac12\int\rho|v|^{10}|\nabla v|^2\dd x+CM^{1/6}Y_{12}^{5/6}.
\]
By H\"older's inequality and \eqref{eq:tm-field12}, the electric work is at most $CY_{12}^{11/12}$. Thus $Y_{12}'\le C(1+Y_{12})$ and
\begin{equation}\label{eq:tm-Y12}
 \sup_{0\le t\le T}Y_{12}(t)\le C_L.
\end{equation}

Next test \eqref{eq:heat-mass} by $\rho^3$. Integration by parts gives
\begin{equation}\label{eq:rho4-moment}
\begin{aligned}
 \frac14\frac d{dt}\int\rho^4\dd x
 +3\int\rho^2|\nabla\rho|^2\dd x
 &=3\int\rho^3v\cdot\nabla\rho\dd x\\
 &\le\tfrac32\int\rho^2|\nabla\rho|^2\dd x
     +\tfrac32\int\rho^4|v|^2\dd x.
\end{aligned}
\end{equation}
Using \eqref{eq:tm-Y12} and then the first bound of \eqref{eq:inputs},
\begin{align*}
 \int\rho^4|v|^2\dd x
 &=\int(\rho|v|^{12})^{1/6}\rho^{23/6}\dd x\\
 &\le Y_{12}^{1/6}\left(\int\rho^{23/5}\dd x\right)^{5/6}\\
 &\le C_L\left(\|\rho\|_\infty^{8/5}
                   \int\rho^3\dd x\right)^{5/6}
 \le C_L\|\rho\|_\infty^{4/3}.
\end{align*}
Integrating \eqref{eq:rho4-moment} and using \eqref{eq:inputs} gives
\begin{align*}
 \|\rho(t)\|_4^4
 &\le\|\rho_0\|_4^4+6C_L\int_0^t\|\rho\|_\infty^{4/3}\dd\tau\\
 &\le\|\rho_0\|_4^4+6C_L
   \Bigl(\int_0^t\|\rho\|_\infty^2\dd\tau\Bigr)^{2/3}t^{1/3}
 \le C_L ,
\end{align*}
hence $\sup_{t\le T}\|\rho(t)\|_4\le C_L$.
H\"older's inequality with $1/12+(11/12)/4=5/16$ yields
\begin{align*}
 \|\rho v\|_{16/5}
 &=\|(\rho^{1/12}v)\rho^{11/12}\|_{16/5}\\
 &\le\|\rho^{1/12}v\|_{12}\|\rho^{11/12}\|_{48/11}
 =Y_{12}^{1/12}\|\rho\|_4^{11/12}\le C_L,\\
 \|\rho(t)\|_\infty
 &\le\|\rho_0\|_\infty+C_L\int_0^t(1+(t-\tau)^{-31/32})\dd\tau\\
 &\le\|\rho_0\|_\infty+C_L(L+32L^{1/32}).
\end{align*}
This proves the required upper bound.
\end{proof}

The moment tests above are justified on each compact interval of positive strong existence by spatial regularization and the Sobolev chain rule. Positivity on that interval justifies the computation. The constants depend on $c_0^{-1}$ through $v_0$.

\subsection{The fourth velocity moment and the relative fifth density moment}\label{sec:upper}
The velocity estimate used here applies on either domain, a fact we will
also use in the lower bound. For the density estimate in this subsection
we take $\Omega=\R^3$ and $\Phi(z)=z^5-b^5-5b^4(z-b)$ for $z\ge0$.
Let $p\ge4$ be real. On any compact interval of strong
existence the positive bounds in \eqref{eq:compact-positive} and the
regularity from Appendix~\ref{app:local} give $v\in C_SH^2$ and
$v_t\in C_SL^2$; hence $Y_p$ is absolutely continuous. Moreover
$|v|^{p-2}v\in H^1$, since its derivative is bounded by
$(p-1)|v|^{p-2}|\nabla v|$ and $v\in L^\infty$, $\nabla v\in L^2$.
For the integrations by parts on $\R^3$, write $A=|v|^{p-2}v$.
Since $p\ge4$ and $v\in L^2\cap L^\infty$, we have
$\|A\|_1\le\|v\|_\infty^{p-3}\|v\|_2^2$ and $A\in L^2$.
The pressure flux $\rho^\gamma A$ and the transport flux $\rho u|v|^p$
are in $L^1$, since $\rho$ and $u$ are bounded and $|v|^p\in L^1$.
The diffusion flux is also integrable, by
$\|\rho(\nabla v)^{\mathsf T}A\|_1\le\|\rho\|_\infty\|\nabla v\|_2\|A\|_2$.
With the cutoffs $\chi_N$ from Lemma~\ref{lem:quantum}, the boundary term
for each flux $F$ is bounded by $N^{-1}\|\nabla\chi\|_\infty\|F\|_1$.
Moreover, $\Div A\in L^1$ since $|v|^{p-2}$ and $\nabla v$ are in $L^2$.
Applying the same cutoff to $A$ gives $\int\Div A\dd x=0$, which allows us
to replace $\Pi$ by $\rho^\gamma$ in the pressure term.
Testing \eqref{eq:effective-velocity} with
$|v|^{p-2}v$ and using \eqref{eq:mass} gives, for almost every time,
\begin{align}
 \tfrac1pY_p'
 &+\int\rho|v|^{p-2}\bigl(|\nabla v|^2+(p-2)|\nabla |v||^2\bigr)\dd x
 \notag\\
 &=\int\rho^\gamma\bigl(|v|^{p-2}\Div v
    +(p-2)|v|^{p-3}v\cdot\nabla|v|\bigr)\dd x
   -\int\rho |v|^{p-2}v\cdot E\dd x.
 \label{eq:general-moment-test}
\end{align}
Here $\nabla|v|=0$ almost everywhere on $\{v=0\}$. Estimate the two pressure
terms separately, using
$|\Div v|\le\sqrt3\,|\nabla v|$ for the first and
$|v\cdot\nabla|v||\le |v|\,|\nabla|v||$ for the second. Young's
inequality absorbs half of each corresponding dissipation term and
leaves coefficients $3/2$ and $(p-2)/2$, respectively. Consequently
\begin{equation}\label{eq:moment-work}
 \frac1pY_p'\le\frac{p+1}{2}\int\rho^{2\gamma-1}|v|^{p-2}\dd x
       +\left|\int\rho |v|^{p-2}v\cdot E\dd x\right|.
\end{equation}

For $p=4$, the pressure term in \eqref{eq:moment-work} carries the
weight $\rho^{2\gamma-1}$. At $\gamma=3$ it has precisely the same
weight as the mixed term from testing the continuity equation by
$\Phi'(\rho)$. The identity
\begin{equation}\label{eq:joint-square}
 \rho^2\bigl(\rho^5+\rho|v|^4\bigr)-2\rho^5|v|^2
 =\rho^3\bigl(\rho^2-|v|^2\bigr)^2
\end{equation}
gives a pointwise upper bound for $\rho^5|v|^2$. Integrating over
$\{\rho>\max(1,2b)\}$ and using $\rho^5\le(16/13)\Phi(\rho)$ there gives
\eqref{eq:high-moments} below. On the complementary set, a constant multiple
of $\rho|v|^2$ bounds the integrand. The coefficient $R^2$ in the high-density
estimate is integrable in time by \eqref{eq:density-time-budget}.

\begin{proposition}\label{prop:upper}
Let $b,\ell>0$ and $1\le\gamma\le3$. Assume the initial data satisfy
\eqref{eq:data}, and let a positive strong solution of the whole-space problem
of Theorem~\ref{thm:cauchy} be given on $[0,S]$. Let $I=[a,a+h]\subset[0,S]$,
where $0<h\le T<\infty$, and assume $\|\rho(a)\|_\infty\le R_{\rm in}$,
$\|v(a)\|_\infty\le W_{\rm in}$. Then
\begin{equation}\label{eq:finite-upper}
 \sup_{t\in I}\|\rho(t)\|_\infty\le M_T
\end{equation}
for a constant $M_T$ depending only on
$T,b,\ell,\gamma,\budget,R_{\rm in}$, $W_{\rm in}$ and absolute constants.
\end{proposition}
\begin{proof}
Write $J_5(t)=\int\Phi(\rho(t))\dd x$. The factorization
\[
 \Phi(z)=(z-b)^2(z^3+2bz^2+3b^2z+4b^3),\qquad z\ge0,
\]
shows that $J_5$ is finite on $I$. Testing the continuity equation in
\eqref{eq:heat-mass}--\eqref{eq:heat-momentum} by $\Phi'(\rho)\in H^1$ gives
\begin{align}
 J_5'+20\int\rho^3|\nabla\rho|^2\dd x
 &=20\int\rho^4v\cdot\nabla\rho\dd x \notag\\
 &\le 10\int\rho^3|\nabla\rho|^2\dd x+10\int\rho^5|v|^2\dd x.
 \label{eq:fifth-moment-absorbed}
\end{align}
For $p=4$, the electric work in \eqref{eq:moment-work} obeys
\begin{equation}\label{eq:fourth-field}
 4\left|\int\rho|v|^2v\cdot E\dd x\right|
 \le4Y_4^{3/4}(R\|E\|_4^4)^{1/4}
 \le3Y_4+R\|E\|_4^4.
\end{equation}
The last field norm is bounded in terms of $b,\gamma,\ell,\budget$ by
\eqref{eq:charge-budget} and \eqref{eq:field-four}.

To bound the mixed terms, set $K_b=\max\{1,2b\}$. On $\{\rho\le K_b\}$,
\[
 (\rho^{2\gamma-1}+\rho^5)|v|^2\le2K_b^4\rho|v|^2
\]
because $0\le2\gamma-2\le4$. On $\{\rho>K_b\}$ we have
$\rho^{2\gamma-1}\le\rho^5$ and $\rho^5\le(16/13)\Phi(\rho)$.
Indeed, $1-5y^4+4y^5\ge13/16$ for $0\le y=b/\rho\le1/2$.
By \eqref{eq:joint-square},
\begin{equation}\label{eq:high-moments}
 \int_{\{\rho>K_b\}}(\rho^{2\gamma-1}+\rho^5)|v|^2\dd x
 \le R^2\bigl((16/13)J_5+Y_4\bigr).
\end{equation}
This bound covers the case $\gamma=3$, where the two density powers coincide. Combining \eqref{eq:moment-work} at $p=4$ with
\eqref{eq:fifth-moment-absorbed}--\eqref{eq:high-moments} and
$Y_2\le2\budget$, we obtain
\begin{equation}\label{eq:joint-gronwall}
 Z'\le C(1+R^2)Z,\qquad Z=1+Y_4+J_5 ,
\end{equation}
where $C=C(b,\gamma,\ell,\budget)$. At $t=a$, the factorization of $\Phi$
and $b\le R_{\rm in}$ give
$\Phi(\rho)\le10R_{\rm in}^3(\rho-b)^2$, hence
$J_5(a)\le10R_{\rm in}^3\|r(a)\|_2^2$. Also
$Y_4(a)\le W_{\rm in}^2Y_2(a)\le2\budget W_{\rm in}^2$.
Gronwall's inequality, \eqref{eq:charge-budget} and
\eqref{eq:density-time-budget} thus give
\begin{equation}\label{eq:fourth-fifth-bound}
 \sup_I Z\le C(1+R_{\rm in}^3+W_{\rm in}^2)e^{C(1+T)}.
\end{equation}

Set $X=\sup_I R$ and $Y=\sup_I Y_4$. These quantities are finite
by local regularity. In the continuity equation we have
$\|m\|_4^4\le X^3Y_4$, since $R(t)\le X$. The free heat evolution preserves the constant
background, and the gradient kernel from $L^4$ to $L^\infty$ has time
singularity $h^{-7/8}$: by \eqref{eq:heat-kernel},
$\|\partial_jG_{t-\tau}*m\|_\infty\le(C_{\rm heat}/8)(t-\tau)^{-7/8}\|m\|_4$
componentwise, and the factor $8$ from
$\int_a^t(t-\tau)^{-7/8}\dd\tau=8(t-a)^{1/8}$ cancels the factor $8$ in the
definition \eqref{eq:gaussian-constant} of $C_{\rm heat}$, with the
remaining dimensional constants absorbed. Therefore \eqref{eq:heat-integral} gives
\begin{equation}\label{eq:mass-kernel-upper}
 X\le R_{\rm in}+C_{\rm heat} T^{1/8}X^{3/4}Y^{1/4}.
\end{equation}
Young's inequality bounds the last term by
$X/2+(27/32)C_{\rm heat}^4T^{1/2}Y$. Inserting
\eqref{eq:fourth-fifth-bound} proves \eqref{eq:finite-upper}.
\end{proof}

\section{A positive lower bound on the density in both geometries}\label{sec:lower}
Let $I=[a,a+h]$ be an interval of positive strong existence, with
$0<h\le T$. Write $M_T$, $W_{\rm in}$ and $K_2\ge1$ for bounds on
$\sup_I R$, $\|v(a)\|_\infty$ and $\sup_I Y_2$, respectively.
On the torus we may take $a=0$ and use
Proposition~\ref{prop:upper-poly} for $1<\gamma\le3$, or
Proposition~\ref{prop:upper-iso} for $\gamma=1$. Conservation of momentum
gives
$Y_2=\int_\T\rho|v-U|^2\dd x+M|U|^2$, so that
$K_2=1+2E_v(0)+M|U|^2$ suffices. In the whole space
Proposition~\ref{prop:upper} provides $M_T$ on any interval,
and \eqref{eq:basic-bounds} permits $K_2=2\budget$.
The corresponding field bound is
\begin{equation}\label{eq:field-finite-bound}
 \sup_I\|E\|_\infty\le F_T,\qquad
 F_T=\begin{cases}
 C_\ell(M_T+M), &\Omega=\T,\\
 \ell^{-2}\bigl(M_T+b+(4\pi)^{-1/2}K_r\bigr),&\Omega=\R^3,
 \end{cases}
\end{equation}
by Lemma~\ref{lem:poisson-estimates} and \eqref{eq:newton-field},
where $K_r$ is the right-hand side of \eqref{eq:charge-budget}.
With the field and the density bounded above, the pressure and electric
work in the velocity equation control the growth of all high moments.

\begin{lemma}\label{lem:large-moments}
Let $I=[a,a+h]\subset[0,S]$ with $0<h\le T$, and suppose that
\[
 \sup_{t\in I}\|\rho(t)\|_\infty\le M_T,\qquad
 \|v(a)\|_\infty\le W_{\rm in},\qquad
 \sup_{t\in I}Y_2(t)\le K_2,\qquad
 \sup_{t\in I}\|E(t)\|_\infty\le F_T ,
\]
where $K_2\ge1$. Then for every real $p\ge6$,
\begin{equation}\label{eq:sqrt-p}
 \sup_{t\in I}Y_p(t)^{1/p}\le A_T\sqrt p,\qquad
 A_T=4K_2^{1/6}
 \bigl(1+W_{\rm in}+M_T^{\gamma-1}\sqrt T+F_TT\bigr).
\end{equation}
\end{lemma}
\begin{proof}
Since $\rho^{2\gamma-1}\le M_T^{2(\gamma-1)}\rho$,
\eqref{eq:moment-work} and \eqref{eq:field-finite-bound} give
\begin{equation}\label{eq:moment-ode}
 \tfrac1pY_p'\le\tfrac{p+1}{2}M_T^{2(\gamma-1)}Y_{p-2}
                       +F_TY_{p-1}.
\end{equation}
Interpolation between $Y_2$ and $Y_p$ works on both domains, including the
case of infinite background mass. Applying H\"older's inequality with measure $\rho\dd x$ gives
\begin{equation}\label{eq:moment-interpolation}
 Y_{p-j}\le Y_2^{j/(p-2)}Y_p^{1-j/(p-2)}
 \le K_2\left(\frac{Y_p}{K_2}\right)^{1-j/(p-2)},
 \qquad j=1,2.
\end{equation}
Set $\Theta=(1+Y_p/K_2)^{2/(p-2)}$. Substitution in
\eqref{eq:moment-ode} gives
\begin{align}
 \Theta'&\le
 \frac{p(p+1)}{p-2}M_T^{2(\gamma-1)}
 +\frac{2p}{p-2}F_T\sqrt{\Theta}\notag\\
 &\le 2pM_T^{2(\gamma-1)}+3F_T\sqrt{\Theta}.
 \label{eq:moment-transformed}
\end{align}
At time $a$ we have $Y_p(a)\le W_{\rm in}^{p-2}Y_2(a)$ and hence
$\Theta(a)\le 2\max\{1,W_{\rm in}^2\}$. For
$X=(\sup_I\Theta)^{1/2}$, integration and the quadratic formula give
\[
 X^2\le \Theta(a)+2pM_T^{2(\gamma-1)}T+3F_TTX,\qquad
 X\le 3F_TT+\sqrt{\Theta(a)+2pM_T^{2(\gamma-1)}T}.
\]
Since $X\ge1$ and $K_2\ge1$, we have
$\sup_I Y_p^{1/p}\le K_2^{1/6}X$, which gives
\eqref{eq:sqrt-p}.
\end{proof}

To convert these weighted moments into a lower density bound, we choose
their order using the minimum density on a short interval. The factor
$\sqrt p$ in \eqref{eq:sqrt-p} then permits a comparison step of
length proportional to $1/p$.

\begin{proposition}\label{prop:lower}
Let the hypotheses of Lemma~\ref{lem:large-moments} hold, and suppose in
addition that $\rho(x,a)\ge c_{\rm in}>0$. Then
\begin{equation}\label{eq:finite-lower}
 \inf_{(x,t)\in\Omega\times I}\rho(x,t)\ge c_T>0 ,
\end{equation}
where $c_T$ depends only on $c_{\rm in},T,A_T$ and dimensional constants.
\end{proposition}
\begin{proof}
On a closed subinterval $I'=[t_0,t_0+h']\subset I$ of positive length define
\[
 \kappa=\max\{e,\sup_{I'}\|\rho^{-1}\|_\infty\},\qquad p=6+\log\kappa.
\]
By \eqref{eq:compact-positive}, $\kappa$ is finite. Fix this $p$
throughout $I'$. Removing the density weight in \eqref{eq:sqrt-p} costs a
factor $\kappa^{1/p}\le e$, since $\rho\ge\kappa^{-1}$ on $I'$; hence
\begin{equation}\label{eq:inverse-moment}
 \sup_{I'}\|v\|_p\le eA_T\sqrt p.
\end{equation}

On $J=[t_1,t_1+d]\subset I'$ consider the comparison equation
\begin{equation}\label{eq:comparison-heat}
 w_t-\Delta w=-\Div(vw),\qquad w(t_1)=1.
\end{equation}
Its mild form is $w=1+\mathcal Lw$, where
\[
 (\mathcal L\chi)(t)=
 -\sum_j\int_{t_1}^t\partial_jG_{t-\tau}*(v_j\chi)(\tau)\dd\tau.
\]
For $d\le1$, \eqref{eq:heat-kernel} and its periodic counterpart give,
uniformly in $p$,
\begin{align}
\|\mathcal L\chi\|_\infty
&\le\sum_j\int_{t_1}^t\|\partial_jG_{t-\tau}*(v_j\chi)\|_\infty\dd\tau
\notag\\
&\le C\int_0^d(d-\tau)^{-\frac12-\frac{3}{2p}}
   \|v\|_p\|\chi\|_\infty\dd\tau
\le C_HeA_T\sqrt p\,\frac{2p}{p-3}\,d^{\frac{p-3}{2p}}\|\chi\|_\infty ,
\label{eq:L-norm}
\end{align}
where we used \eqref{eq:inverse-moment} and
$\int_0^d(d-\tau)^{-\alpha_p}\dd\tau=d^{1-\alpha_p}/(1-\alpha_p)$ with
$\alpha_p=1/2+3/(2p)$. Set $D_{\rm c}=(1+32C_HeA_T)^4$. If
$d\le1/(D_{\rm c}p)$, then
\[
 \sqrt p\,d^{(p-3)/(2p)}
 \le D_{\rm c}^{-(p-3)/(2p)}p^{3/(2p)}
 \le D_{\rm c}^{-1/4}p^{3/(2p)}
 \le 6^{1/4}D_{\rm c}^{-1/4}<2D_{\rm c}^{-1/4},
\]
because $p\ge6$, $(p-3)/(2p)\ge1/4$, and $(\log p)/p$ decreases on
$[6,\infty)$. Since $2p/(p-3)\le4$, \eqref{eq:L-norm} gives
$\|\mathcal L\|\le1/4$, and the fixed point satisfies
$|w-1|\le1/3$ almost everywhere on $J\times\Omega$.

We compare $\rho$ with $cw$ through the negative part of their difference.
The flux $vw$ belongs to $L^2(J;L^2)$, so \eqref{eq:linear-integrated} at
$k=0$ gives $w-1\in C(J;L^2)\cap L^2(J;H^1)$ with zero initial trace. Set
$c=\inf_x\rho(x,t_1)$ and $q=\rho-cw$. The conserved mean on the torus or the
limit at infinity in the whole space gives $0<c\le b$, and
\[
 q_t-\Delta q=-\Div(vq),\qquad q(t_1)\ge0.
\]
Although $q$ may have a nonzero background $b-c$, the function
$y=q-(b-c)$ belongs to $C(J;L^2)\cap L^2(J;H^1)$, and
$y_t=q_t\in L^2(J;H^{-1})$ because $vq\in L^2(J;L^2)$. The map
$z\mapsto(z+b-c)_-$ is Lipschitz and vanishes at zero, where
$z_-=\max\{-z,0\}$; hence $q_-\in C(J;L^2)\cap L^2(J;H^1)$. Testing the
equation for $q$ with $-q_-$ gives
\[
 \tfrac12\frac d{dt}\|q_-\|_2^2+\|\nabla q_-\|_2^2
 =\int vq_-\cdot\nabla q_-\dd x
 \le\tfrac12\|\nabla q_-\|_2^2+\tfrac12\|v\|_\infty^2\|q_-\|_2^2.
\]
For the time chain rule, first mollify $y$ in time on compact subintervals
of the interior of $J$: the derivatives converge in $L^2(H^{-1})$ and the
shifted negative parts converge in $L^2(H^1)$ by the Sobolev chain rule and
$\nabla q=0$ almost everywhere on $\{q=0\}$. Continuity into $L^2$ of
$q_-$ then gives the endpoint terms. Since $q_-(t_1)=0$ and
$v\in C(J;H^2)$, Gronwall's inequality gives $q_-=0$. This comparison uses
only the finiteness of $\|v\|_{L^\infty(J\times\Omega)}$. We obtain
$\rho\ge cw\ge2c/3$ almost everywhere, and continuity of $\rho$ yields
\begin{equation}\label{eq:comparison-step}
 \inf_{\Omega\times J}\rho\ge\tfrac23\inf_x\rho(x,t_1)
 \qquad\text{whenever } d\le1/(D_{\rm c}p).
\end{equation}

Partition $I'$ into $n=\lceil D_{\rm c}ph'\rceil$ equal pieces. Iterating
\eqref{eq:comparison-step} gives
\[
 \inf_{\Omega\times I'}\rho\ge2^{-n}\inf_{x\in\Omega}\rho(x,t_0).
\]
For $\varphi(t)=\log\max\{e,\|\rho(t)^{-1}\|_\infty\}$ it follows that
\begin{equation}\label{eq:inverse-absorption}
 \log\kappa\le\varphi(t_0)+\log2+D_{\rm c}h'\log2\,(6+\log\kappa).
\end{equation}
If $h'\le h_0=(2D_{\rm c}\log2)^{-1}$, absorption gives
$\sup_{I'}\varphi=\log\kappa\le2\varphi(t_0)+6+2\log2$.
The entire interval $I$ is covered by at most $N_0=\lceil T/h_0\rceil$
successive intervals of this length or less. Since
$\varphi(a)\le\log\max\{e,c_{\rm in}^{-1}\}$, \eqref{eq:finite-lower}
holds with
\[
 c_T=\exp\!\left[-2^{N_0}\bigl(\log\max\{e,c_{\rm in}^{-1}\}+2\log2+6\bigr)\right].
\]
\end{proof}

\section{Higher-order regularity and global continuation}\label{sec:continuation}
The eighth velocity moment controls the convective flux in $L^4$.
The heat equation for $m$ then bounds the remaining pointwise
coefficient in the $H^2$ estimates for $\mathcal U$.

\begin{proposition}\label{prop:global}
Under the data assumptions of Theorem~\ref{thm:periodic} or
Theorem~\ref{thm:cauchy} on the corresponding domain, there is a unique
global positive strong solution in the sense of
Definition~\ref{def:positive-strong-solution}. For each finite $T>0$, it
satisfies \eqref{eq:strong-regularity} or \eqref{eq:regularity}, respectively, and
\[
 \sup_{0\le t\le T}
   \bigl(\|r(t)\|_{H^3}+\|u(t)\|_{H^2}+\|v(t)\|_{H^2}\bigr)\le C_T.
\]
Here $C_T$ depends only on $T$ and the parameters and initial data bounds
listed in the corresponding theorem.
\end{proposition}
\begin{proof}
Proposition~\ref{prop:local-existence} applies to the data by
\eqref{eq:initial-effective-data}. Let $T^*$ be the maximal lifespan of
Proposition~\ref{prop:local-existence}.
Fix $T>0$ and $S<\min\{T,T^*\}$. On $[0,S]$ the upper bounds from
Section~\ref{sec:upper-both}, Lemma~\ref{lem:large-moments} and
Proposition~\ref{prop:lower} give
\begin{equation}\label{eq:finite-density-bounds}
 c_T\le\rho\le M_T,\qquad \|E\|_\infty\le F_T,\qquad
 Y_8^{1/8}\le\sqrt8\,A_T.
\end{equation}
The constants are independent of $S$. In particular
$\|v_0\|_\infty\le S_2(\|u_0\|_{H^2}+c_0^{-1}\|r_0\|_{H^3})$,
and the initial energy bounds depend only on the stated data.
We have
\[
 \|m\otimes m/\rho\|_4^4\le M_T^3Y_8,\qquad
 \|\Pi\|_\infty\le M_T^\gamma+b^\gamma,\qquad
 \|\rho E\|_\infty\le M_TF_T.
\]
For times $0<d\le1$, the differentiated kernel in
\eqref{eq:heat-integral} has time singularity $d^{-7/8}$ from $L^4$ to
$L^\infty$ and $d^{-1/2}$ on $L^\infty$. For $d\ge1$ both operator
norms stay bounded, by the semigroup property and contraction on
$L^\infty$. Since $m_0=\rho_0u_0+\nabla r_0$ is bounded, we obtain
\begin{equation}\label{eq:momentum-sup-bound}
  \sup_{0\le t\le S}\|m(t)\|_\infty\le C_T.
\end{equation}
Together with \eqref{eq:finite-density-bounds} and
\eqref{eq:momentum-sup-bound} this places $\mathcal U$
in a bounded subset of $\{b+r\ge c_T\}$. Along the segments from zero
to each value of $\mathcal U$, one has $b+\theta r\ge\min\{b,c_T\}$.
Each pointwise component $F$ of the flux and source in \eqref{eq:heat-maps}
vanishes at zero and has bounded first and second derivatives on these
segments. Thus
\[
  \|F(\mathcal U)\|_{H^2}
  \le C_T\bigl(\|\mathcal U\|_{H^2}+\|\nabla\mathcal U\|_4^2\bigr).
\]
Integration by parts in $\int |\nabla\mathcal U|^4\dd x$, initially for
smooth periodic or compactly supported functions, gives
\[
\int|\nabla\mathcal U|^4\dd x
=-\int(\nabla\mathcal U\cdot\mathcal U)\cdot\nabla|\nabla\mathcal U|^2\dd x
 -\int(\mathcal U\cdot\Delta\mathcal U)|\nabla\mathcal U|^2\dd x,
\]
and hence
\[
\|\nabla\mathcal U\|_4^4
\le(2+\sqrt3)\|\mathcal U\|_\infty\|\nabla\mathcal U\|_4^2
\|\nabla^2\mathcal U\|_2,
\]
where we used $|\Delta\mathcal U|\le\sqrt3|\nabla^2\mathcal U|$ and
$|\nabla|\nabla\mathcal U|^2|\le2|\nabla\mathcal U||\nabla^2\mathcal U|$.
Dividing by $\|\nabla\mathcal U\|_4^2$ (the case
$\|\nabla\mathcal U\|_4=0$ is trivial) and passing to $H^2$ functions by
density gives
\begin{equation}\label{eq:gradient-product}
 \|\nabla\mathcal U\|_4^2\le(2+\sqrt3)\|\mathcal U\|_\infty
                                    \|\nabla^2\mathcal U\|_2.
\end{equation}
The remaining source component $Qm$ is bounded on $H^2$. Since
$\mathcal U\in C([0,S];H^2)$ by local existence,
$\|\mathcal U\|_\infty\le C_T$ and $\|\nabla^2\mathcal U\|_2\le
\|\mathcal U\|_{H^2}$; by \eqref{eq:gradient-product} the flux bound above
is therefore at most $C_T\|\mathcal U\|_{H^2}$. The source obeys
$\|\mathcal S(\mathcal U)\|_{H^2}\le C_T(1+\|\mathcal U\|_{H^2})$, by the
pointwise bounds in \eqref{eq:finite-density-bounds}, the tame $H^2$
product estimate and \eqref{eq:gauss-symbol}. Applying
\eqref{eq:linear-energy} with $k=2$ gives
\begin{equation}\label{eq:continuation-ode}
 \frac d{dt}\|\mathcal U\|_{H^2}^2+\|\nabla\mathcal U\|_{H^2}^2
   \le C_T\bigl(1+\|\mathcal U\|_{H^2}^2\bigr).
\end{equation}
Gronwall's inequality bounds $\mathcal U$ in
$C([0,S];H^2)\cap L^2(0,S;H^3)$ by the stated data and $T$.
Applying \eqref{eq:linear-integrated} at $k=3$ to the continuity equation
with flux $m$ bounds $r$ in $C([0,S];H^3)\cap L^2(0,S;H^4)$.
Composition on $[c_T,M_T]$ bounds $1/\rho-1/b$ in $C([0,S];H^3)$,
so the identities $v=m/\rho$ and $u=(m-\nabla r)/\rho$ give the
claimed quantitative $H^2$ bounds.

Suppose that $T^*<\infty$. Fixing $S<T^*$ and letting $S\uparrow T^*$ in
the preceding estimates gives a common $H^2$ bound for $\mathcal U$ and the
bound $c_T\le\rho\le M_T$ on $[0,T^*)$. Every trace $\mathcal U(a)$, $a<T^*$,
satisfies the field constraints and has this $H^2$ bound and lower density
bound; on the torus the traces also preserve the zero spatial means of $r$
and $E$. Proposition~\ref{prop:local-existence} therefore assigns a common
positive lifespan $h_{\rm loc}$ to all these traces. Restarting at a time
$a>T^*-h_{\rm loc}/2$ extends the solution past $T^*$, and uniqueness on the
overlap gives a single solution. Hence $T^*=\infty$. The regularity
assertions follow by applying the local result on every compact interval.
For uniqueness, let two positive strong solutions on $[0,T]$ in the sense of
Definition~\ref{def:positive-strong-solution} share the same data. By the
converse implication at the end of the proof of
Proposition~\ref{prop:local-existence}, both are
positive $C([0,T];H^2)$ mild solutions of
\eqref{eq:heat-mass}--\eqref{eq:heat-momentum} and \eqref{eq:field-heat}.
Proposition~\ref{prop:local-existence} then gives a common lifespan at every
common trace, and iterating the argument shows that the two solutions
coincide on $[0,T]$.
\end{proof}

\section{Exponential relaxation on the torus}\label{sec:relaxation}
For the global solution constructed in Section~\ref{sec:continuation},
we now derive uniform density bounds and exponential convergence to
$(M,U,0)$. The argument applies throughout $1\le\gamma\le3$ and uses
the energy identities of Lemma~\ref{lem:energy} together with the
regularity on finite time intervals. The density bounds uniform in
time are obtained as part of the conclusion.

\subsection{Approach to equilibrium}\label{sec:equilibrium-entry}
With $\mathcal I$ as in Lemma~\ref{lem:quantum}, define
\[
 G(t)=\int\rho|\nabla v|^2\dd x,
 \qquad B_0=E_u(0)+3E_v(0).
\]
Adding the bounds for $\mathcal I$ and $G$ in \eqref{eq:Iintegral} and \eqref{eq:Ev} and letting the terminal time tend to infinity gives
\begin{equation}\label{eq:tail-budget}
 \int_0^\infty (\mathcal I+G)\dd t\le B_0.
\end{equation}

The following consequence of \eqref{eq:tail-budget} and the
conservation laws provides the connection to the convergence estimates.

\begin{lemma}\label{lem:good-times}
For the global solution above and every $\varepsilon>0$, there is a
time $t_*\in(0,L_\varepsilon)$ such that
\begin{gather*}
 \|\rho(t_*)-M\|_{H^2}+\|v(t_*)-U\|_{H^1}\le\varepsilon,\\
 M/2\le\rho(x,t_*)\le3M/2\qquad(x\in\T).
\end{gather*}
The time bound may be chosen as
\begin{equation}\label{eq:tail-entry-time}
 L_\varepsilon=1+2B_0/\eta,
\end{equation}
where $\eta>0$ depends only on $M$, $\varepsilon$, and the unit torus.
Moreover, there is a sequence $t_n\to\infty$ such that
\[
 \|\rho(t_n)-M\|_{H^2}+\|v(t_n)-U\|_{H^1}\to0.
\]
\end{lemma}
\begin{proof}
The nonzero frequencies and \eqref{eq:quantum} give
$\|\sqrt\rho-\overline{\sqrt\rho}\|_{H^2}\le C \mathcal I(t)^{1/2}$.
Mass conservation implies
\[
 M=(\overline{\sqrt\rho})^2+\|\sqrt\rho-\overline{\sqrt\rho}\|_2^2,\qquad
 |\overline{\sqrt\rho}-\sqrt M|\le C M^{-1/2}\mathcal I(t).
\]
If $\mathcal I(t)\le1$, the algebra property of $H^2$ applied to
$\rho-M=(\sqrt\rho-\sqrt M)(\sqrt\rho+\sqrt M)$ yields
\begin{equation}\label{eq:tail-good-density}
 \|\rho(t)-M\|_{H^2}\le C_M \mathcal I(t)^{1/2}
 \qquad\text{whenever }\mathcal I(t)\le1.
\end{equation}
Thus there is $\eta_0=\eta_0(M)\in(0,1]$ such that $\mathcal I(t)\le\eta_0$
implies $M/2\le\rho(t)\le3M/2$.
Put $w=v-U$ and $\bar w=\int w\dd x$. Since
$\int\rho w\dd x=0$ by Lemma~\ref{lem:energy} and
$\int(\rho-M)\dd x=0$, Poincar\'e's inequality yields
\begin{equation}\label{eq:tail-mean}
 M\bar w=-\int(\rho-M)(w-\bar w)\dd x,
 \qquad
 |\bar w|\le C M^{-1}\|\rho-M\|_2\|\nabla w\|_2.
\end{equation}
When $\mathcal I(t)\le\eta_0$, the preceding density bounds imply
$\|\nabla w\|_2^2\le2G(t)/M$. By \eqref{eq:tail-good-density}, \eqref{eq:tail-mean} and Poincar\'e's inequality,
\[
 \|\rho-M\|_{H^2}+\|w\|_{H^1}\le C_M(\mathcal I(t)+G(t))^{1/2}.
\]
Choose $0<\eta\le\eta_0$ so that $C_M\sqrt\eta\le\varepsilon$. By \eqref{eq:tail-budget}, the part of $(0,L_\varepsilon)$ on which $\mathcal I+G>\eta$ has measure at most $B_0/\eta<L_\varepsilon$. Select $t_*$ in its complement. This proves both estimates at $t_*$.
If $B_0=0$, then $\mathcal I+G=0$ almost everywhere, and the same argument
applies with any $t_*\in(0,L_\varepsilon)$.
To obtain the sequence, choose
$t_n\in(n,n+1)$ with
$\mathcal I(t_n)+G(t_n)\le2\int_n^{n+1}(\mathcal I+G)\dd t$. The right-hand side tends to
zero by integrability.
\end{proof}

\subsection{Exponential stability}\label{sec:exponential-stability}
We next establish quantitative evolution estimates near $(M,U)$.
Lemma~\ref{lem:good-times} supplies a time at which these estimates
apply to the global solution.

\begin{proposition}\label{prop:tail-stability}
On the unit torus $\T$, fix $M>0$, $1\le\gamma\le3$ and $\ell>0$. There are
$\varepsilon_*>0$, $C>0$ and $\kappa>0$, depending only on these
parameters and the torus, with the following property. For any $U\in\R^3$,
suppose a positive strong solution has mass $M$, momentum $MU$, and at
some time $t_*$ in its interval of existence
satisfies
\[
 \|\rho(t_*)-M\|_{H^2}+\|v(t_*)-U\|_{H^1}\le\varepsilon_*.
\]
Throughout its subsequent interval of strong existence,
\begin{align}
 M/2&\le\rho(x,t)\le3M/2,\label{eq:tail-density}\\
 \|\rho(t)-M\|_{H^2}^2+\|v(t)-U\|_{H^1}^2
 &\le C\varepsilon_*^2 e^{-\kappa(t-t_*)},\label{eq:tail-low}\\
 \|\rho(t)-M\|_{H^3}^2+\|v(t)-U\|_{H^2}^2
 &\le C A_* e^{-\kappa(t-t_*)},\label{eq:tail-high}
\end{align}
where
$A_*=\|\rho(t_*)-M\|_{H^3}^2+\|v(t_*)-U\|_{H^2}^2$.
If the solution is global, then
\begin{equation}\label{eq:tail-high-integral}
 \int_{t_*}^\infty
  \bigl(\|\rho-M\|_{H^4}^2+\|v-U\|_{H^3}^2\bigr)\dd t
 \le C A_*.
\end{equation}
\end{proposition}

\begin{proof}
Set $q=\rho-M$ and $w=v-U$, and change variables by
$y=x-U(t-t_*)$, that is, consider
$\tilde q(y,t)=q(y+U(t-t_*),t)$ and
$\tilde w(y,t)=w(y+U(t-t_*),t)$; denote the new functions again by $q$ and
$w$. Then $\partial_t$ becomes $\partial_t-U\cdot\nabla$, and the terms
containing $U$ cancel. The effective equations become
\begin{align}
 q_t-\Delta q+M\Div w&=-\Div(qw),\label{eq:tail-q}\\
 w_t-\Delta w+MB\nabla q&=N,\label{eq:tail-w}\\
 B&=M^{-1}\bigl(\gamma M^{\gamma-2}\operatorname{Id}+\ell^{-2}(-\Delta)^{-1}\bigr),\label{eq:tail-B}\\
 N&=-w\cdot\nabla w+2(\nabla g\cdot\nabla)w-\gamma\bigl((M+q)^{\gamma-2}-M^{\gamma-2}\bigr)\nabla q,\label{eq:tail-N}\\
 g&=\log(M+q).\notag
\end{align}
Here $\operatorname{Id}$ is the identity operator and $g=\log\rho$.
The inverse Laplacian acts on mean-zero periodic functions, and $B$
acts componentwise when applied to a vector field. The constraints are
\begin{equation}\label{eq:tail-constraints}
 \int q\dd x=0,\qquad \int(M+q)w\dd x=0.
\end{equation}
The first follows from mass conservation. For the second, we use
$\int_\T\rho v=\int_\T\rho u=MU$ from Lemma~\ref{lem:energy}, since
$\int_\T\rho\nabla\log\rho\dd x=\int_\T\nabla\rho\dd x=0$.
This change of variables preserves all the norms used below, and the constant
velocity $U$ no longer appears explicitly in the equations.

Use the multipliers $J$ and $\Lambda$ defined in Section~\ref{sec:effective}.
On mean-zero functions the operator $B$ from \eqref{eq:tail-B} has symbol
\[
 \beta(\xi)=\frac{\gamma M^{\gamma-2}+\ell^{-2}|\xi|^{-2}}{M},
 \qquad \xi\in2\pi\mathbb Z^3\setminus\{0\}.
\]
In particular,
\[
 0<\beta_-:=\gamma M^{\gamma-3}\le \beta(\xi)
 \le\gamma M^{\gamma-3}+\frac1{4\pi^2M\ell^2}=:\beta_+<\infty.
\]
Fix $\lambda=[16(1+M^2)]^{-1}$, so that $\lambda M^2\le1/16$.
Define, for $j=1,2$,
\begin{align}
 \mathcal E_j&=\frac12\left(
   \langle BJ^jq,J^jq\rangle_{L^2}+\|J^jw\|_2^2
   +\lambda\|\Lambda^{j+1}q\|_2^2\right),\label{eq:tail-E}\\
 \mathcal D_j&=\langle B\nabla J^jq,\nabla J^jq\rangle_{L^2}
   +\|\nabla J^jw\|_2^2
   +\lambda\|\Lambda^{j+2}q\|_2^2.\label{eq:tail-D}
\end{align}
For the Sobolev orders $j=1,2$, set
$X_j=(\|q\|_{H^{j+1}}^2+\|w\|_{H^j}^2)^{1/2}$.
Under the temporary assumption
\begin{equation}\label{eq:tail-bootstrap}
 X_1\le\delta\le1,
\end{equation}
the mean identity \eqref{eq:tail-mean} and Poincar\'e's inequality give $\|w\|_{H^s}\le C_M\|\nabla w\|_{H^{s-1}}$ for $s=1,2,3$. The bounds $\beta_-\le \beta(\xi)\le \beta_+$ and the spectral gap on mean-zero functions therefore imply
\begin{align}
 \mathcal E_j&\asymp X_j^2,\notag\\
 \mathcal D_j&\asymp\|q\|_{H^{j+2}}^2+\|w\|_{H^{j+1}}^2,
 \qquad \mathcal E_j\le C\mathcal D_j.\label{eq:tail-equivalence}
\end{align}
Test the density and velocity equations at order $j$ against $BJ^jq$ and $J^jw$, respectively, after applying $J^j$. Integration by parts cancels the two terms with coefficient $M$. Testing the continuity equation also at order $j+1$ against $\lambda\Lambda^{j+1}q$ yields
\begin{align}
 \mathcal E_j'+\mathcal D_j
 ={}&-\langle BJ^j\Div(qw),J^jq\rangle_{L^2}
       +\langle J^jN,J^jw\rangle_{H^{-1},H^1}\notag\\
 &-\lambda M\langle\Lambda^{j+1}\Div w,\Lambda^{j+1}q\rangle_{H^{-1},H^1}
 -\lambda\langle\Lambda^{j+1}\Div(qw),\Lambda^{j+1}q\rangle_{H^{-1},H^1}.
 \label{eq:tail-identity}
\end{align}
The linear remainder obeys
\begin{align}
 \lambda M\bigl|\langle\Lambda^{j+1}\Div w,
                         \Lambda^{j+1}q\rangle_{H^{-1},H^1}\bigr|
 &=\lambda M\bigl|\langle\Lambda^{j+1}w,
                         \nabla\Lambda^{j+1}q\rangle_{L^2}\bigr|\notag\\
 &\le\lambda M\|\Lambda^{j+1}w\|_2\|\Lambda^{j+2}q\|_2\notag\\
 &\le\frac\lambda8\|\Lambda^{j+2}q\|_2^2
       +2\lambda M^2\|\Lambda^{j+1}w\|_2^2\notag\\
 &\le\tfrac18\lambda\|\Lambda^{j+2}q\|_2^2
       +\tfrac18\|\nabla J^jw\|_2^2
 \le\tfrac14\mathcal D_j.\label{eq:tail-linear}
\end{align}
The last two inequalities follow frequency by frequency from
$\|\Lambda^{j+1}w\|_2\le\|\nabla J^jw\|_2$ and from $2\lambda M^2\le1/8$.
The mean of $w$
does not enter this estimate. The other terms satisfy
\begin{align}
 |\langle BJ^j\Div(qw),J^jq\rangle_{L^2}|
   &=|\langle J^j(qw),B\nabla J^jq\rangle_{L^2}|\notag\\
   &\le\|J^j(qw)\|_2\|B\nabla J^jq\|_2\notag\\
   &\le \beta_+\|qw\|_{H^j}\|q\|_{H^{j+1}},\label{eq:tail-F-low}\\
 |\langle J^jN,J^jw\rangle_{H^{-1},H^1}|
   &=|\langle J^{j-1}N,J^{j+1}w\rangle_{L^2}|\notag\\
   &\le\|N\|_{H^{j-1}}\|w\|_{H^{j+1}},\label{eq:tail-N-dual}\\
 \lambda|\langle\Lambda^{j+1}\Div(qw),\Lambda^{j+1}q\rangle_{H^{-1},H^1}|
   &=\lambda|\langle\Lambda^{j+1}(qw),\nabla\Lambda^{j+1}q\rangle_{L^2}|\notag\\
   &\le\lambda\|\Lambda^{j+1}(qw)\|_2\|\Lambda^{j+2}q\|_2\notag\\
   &\le C\lambda\|qw\|_{H^{j+1}}\|q\|_{H^{j+2}}.
   \label{eq:tail-F-high}
\end{align}
These are $H^{-1}$--$H^1$ pairings at the highest order. On compact intervals $q\in L^2H^4$, $q_t\in L^2H^2$, $w\in L^2H^3$ and $w_t\in L^2H^1$. Finite Fourier projections and the energy-continuity property in Section~\ref{sec:effective} justify the identity between any two times on such intervals.

\emph{First level.}
Choose $\delta$ small enough that \eqref{eq:tail-bootstrap} implies
$|q|\le M/2$. Leibniz' rule and the three-dimensional embeddings give
\begin{align}
 \|qw\|_{H^1}&\le C\|q\|_{H^2}\|w\|_{H^1},&
 \|qw\|_{H^2}&\le C\|q\|_{H^2}\|w\|_{H^2}.\label{eq:tail-products-low}
\end{align}
Here $H^2$ is an algebra and acts by multiplication on $H^1$. Under $X_1\le\delta$ the three terms of $N$ satisfy
\begin{align*}
 \|w\cdot\nabla w\|_2
   &\le\|w\|_6\|\nabla w\|_3
     \le C\|w\|_{H^1}\|w\|_{H^2},\\
 \|(\nabla g\cdot\nabla)w\|_2
   &\le\|\nabla g\|_6\|\nabla w\|_3
     \le C\|q\|_{H^2}\|w\|_{H^2},\\
 \|\gamma((M+q)^{\gamma-2}-M^{\gamma-2})\nabla q\|_2
   &\le C\|q\|_\infty\|\nabla q\|_2.
\end{align*}
Hence
\begin{equation}\label{eq:tail-N-low}
 \|N\|_2\le C\delta(\|w\|_{H^2}+\|q\|_{H^2}).
\end{equation}
We bound the three nonlinear terms in \eqref{eq:tail-identity} at $j=1$ separately. By \eqref{eq:tail-F-low} and \eqref{eq:tail-products-low}, the flux pairing obeys
\[\begin{aligned}
 \bigl|\langle BJ\Div(qw),Jq\rangle_{L^2}\bigr|
 &=\bigl|\langle J(qw),B\nabla Jq\rangle_{L^2}\bigr|\\
 &\le\sqrt{\beta_+}\,\|qw\|_{H^1}\langle B\nabla Jq,\nabla Jq\rangle_{L^2}^{1/2}\\
 &\le C\|q\|_{H^2}\|w\|_{H^1}\mathcal D_1^{1/2}\\
 &\le CX_1^2\mathcal D_1^{1/2},
\end{aligned}\]
since $\langle B\nabla Jq,\nabla Jq\rangle_{L^2}\le\mathcal D_1$. The first two terms of $N$, paired with $J^2w$ as in \eqref{eq:tail-N-dual}, and the pairing in \eqref{eq:tail-F-high} contribute $CX_1\mathcal D_1$: by the first two bounds above and \eqref{eq:tail-products-low},
\[
 \bigl|\langle w\cdot\nabla w,J^2w\rangle_{L^2}\bigr|
 \le\|w\|_6\|\nabla w\|_3\|J^2w\|_2
 \le C\|w\|_{H^1}\|w\|_{H^2}^2\le CX_1\mathcal D_1,
\]
and similarly for $(\nabla g\cdot\nabla)w$, while
\eqref{eq:tail-F-high} gives
$C\lambda\|q\|_{H^2}\|w\|_{H^2}\|q\|_{H^3}\le C\lambda X_1\mathcal D_1$,
where we used $\|w\|_{H^2}^2+\|w\|_{H^2}\|q\|_{H^3}\le C\mathcal D_1$ from
\eqref{eq:tail-equivalence}. On $M/2\le M+q\le3M/2$, the mean value theorem gives
$|\gamma((M+q)^{\gamma-2}-M^{\gamma-2})|\le C_{M,\gamma}|q|$.
Multiplication by $\nabla q$ therefore gives the quadratic pressure remainder in $N$, with
\[
 \bigl|\langle\gamma((M+q)^{\gamma-2}-M^{\gamma-2})\nabla q,J^2w\rangle_{L^2}\bigr|
 \le C\|q\|_\infty\|\nabla q\|_2\|w\|_{H^2}
 \le CX_1^2\mathcal D_1^{1/2},
\]
where $\|q\|_\infty\le S_2\|q\|_{H^2}\le C\delta$, $\|\nabla q\|_2\le X_1$
and $\|w\|_{H^2}\le C\mathcal D_1^{1/2}$. By Young's inequality, each of the
two bounds of order $X_1^2\mathcal D_1^{1/2}$ is at most
$\tfrac1{32}\mathcal D_1+8C^2X_1^4$. Combining these bounds with
\eqref{eq:tail-linear} and shrinking $\delta$ so that $C\delta\le3/16$,
we obtain
\begin{equation}\label{eq:tail-E1}
  \mathcal E_1'+\tfrac12\mathcal D_1\le C X_1^4 .
\end{equation}
Write $c_E X_1^2\le\mathcal E_1\le C_E X_1^2$ and
$\mathcal E_1\le C\mathcal D_1$ by \eqref{eq:tail-equivalence}. Since
$X_1^4\le c_E^{-2}\mathcal E_1^2$, \eqref{eq:tail-E1} gives
\[
  \mathcal E_1'\le-\kappa_1\mathcal E_1+C_1\mathcal E_1^2,
\]
with parameter-dependent $\kappa_1>0$ and $C_1\ge1$ (explicitly,
$\kappa_1=1/(2C)$ and $C_1=C/c_E^2$). Choose $\varepsilon_*$ so that
$C_E\varepsilon_*^2\le\min\{c_E\delta^2/8,\kappa_1/(2C_1)\}$. Under the
temporary assumption $X_1\le\delta$, the logistic bound
\[
  \mathcal E_1(t)
  \le\frac{\kappa_1\mathcal E_1(t_*)e^{-\kappa_1(t-t_*)}}
       {\kappa_1-C_1\mathcal E_1(t_*)\bigl(1-e^{-\kappa_1(t-t_*)}\bigr)}
  \le2\mathcal E_1(t_*)e^{-\kappa_1(t-t_*)}
  \le\frac{c_E\delta^2}{4}
\]
holds, because $C_1\mathcal E_1(t_*)\le C_1C_E\varepsilon_*^2\le\kappa_1/2$
and $2C_E\varepsilon_*^2\le c_E\delta^2/4$. Hence $X_1\le\delta/2$, which
improves the assumed bound. Integrating \eqref{eq:tail-E1} and using this
decay together with $X_2^2\le C\mathcal D_1$, we obtain
\begin{equation}\label{eq:tail-E1-decay}
  \mathcal E_1(t)\le2\mathcal E_1(t_*)e^{-\kappa_1(t-t_*)},
  \qquad
  \int_{t_*}^t X_2(\tau)^2\dd\tau\le C\mathcal E_1(t_*).
\end{equation}
The time integral controls the coefficient in the next energy inequality.

\emph{Second level.}
The high-order product satisfies
\begin{equation}\label{eq:tail-product-high}
 \|qw\|_{H^3}\le
 C\delta\|w\|_{H^3}+C\|w\|_{H^2}\|q\|_{H^3}.
\end{equation}
Indeed, the tame third-order product bound gives $C(\|q\|_\infty\|w\|_{H^3}+\|w\|_\infty\|q\|_{H^3})$. Since $H^2\hookrightarrow L^\infty$, this proves \eqref{eq:tail-product-high}. Tame composition on $[M/2,3M/2]$ gives
\begin{equation}\label{eq:tail-log}
 \|\nabla g\|_{H^2}\le C\|q\|_{H^3}.
\end{equation}
The multiplier bound $H^2\cdot H^1\subset H^1$ applied to $w\cdot\nabla w$ and $(\nabla g\cdot\nabla)w$, together with \eqref{eq:tail-log}, yields their bounds at order one. The pressure coefficient $\gamma((M+q)^{\gamma-2}-M^{\gamma-2})$ has $H^2$ norm at most $C\|q\|_{H^2}$ by the same composition estimate. Consequently
\begin{equation}\label{eq:tail-N-high}
 \|N\|_{H^1}\le C\left(
   \|w\|_{H^2}^2+\|q\|_{H^3}\|w\|_{H^2}
   +\delta\|q\|_{H^2}\right).
\end{equation}
For $j=2$, \eqref{eq:tail-F-low} and \eqref{eq:tail-products-low} give a bound $C\delta\mathcal D_2$ for the first nonlinear term. By \eqref{eq:tail-N-dual}, \eqref{eq:tail-N-high} and \eqref{eq:tail-equivalence}, the term involving $N$ is bounded by
\[
 C X_2^2\mathcal D_2^{1/2}+C\delta\mathcal D_2
 \le \tfrac18\mathcal D_2+C X_2^2\mathcal E_2+C\delta\mathcal D_2.
\]
The final term in \eqref{eq:tail-F-high} satisfies the same bound by \eqref{eq:tail-product-high}. Decreasing $\delta$ so that the sum of the terms $C\delta\mathcal D_2$ is at most $\mathcal D_2/4$, and using \eqref{eq:tail-linear}, gives
\begin{equation}\label{eq:tail-E2}
 \mathcal E_2'+\tfrac14\mathcal D_2\le C_2X_2^2\mathcal E_2.
\end{equation}
We now fix $\delta$ satisfying both levels and choose $\varepsilon_*$
so that $C_E\varepsilon_*^2\le\min\{c_E\delta^2/8,\kappa_1/(2C_1)\}$. Continuity and
\eqref{eq:tail-E1} extend the improved bound $X_1\le\delta/2$
throughout the interval of strong existence. Thus \eqref{eq:tail-E1-decay}
holds there and bounds
$\int_{t_*}^tX_2(\tau)^2\dd\tau$.
By \eqref{eq:tail-equivalence}, choose $\kappa_2>0$, depending only
on the fixed parameters, such that
$\mathcal D_2/4\ge\kappa_2\mathcal E_2$.
Gronwall's inequality now gives
\[
 \mathcal E_2(t)\le\mathcal E_2(t_*)\exp\left(-\kappa_2(t-t_*)
               +C_2\int_{t_*}^t X_2(\tau)^2\dd\tau\right)
 \le C\mathcal E_2(t_*)e^{-\kappa_2(t-t_*)}.
\]
Indeed, $\int_{t_*}^t X_2^2\dd\tau\le C\mathcal E_1(t_*)\le C\varepsilon_*^2$. Integrating \eqref{eq:tail-E2} gives
\[
 \int_{t_*}^t\mathcal D_2\dd\tau\le C(\mathcal E_2(t_*)+
\sup_{[t_*,t]}\mathcal E_2\int_{t_*}^tX_2^2\dd\tau)\le C A_*.
\]
Equivalence \eqref{eq:tail-equivalence} yields \eqref{eq:tail-high-integral} upon taking $t\to\infty$ for a global solution. The choice $\kappa=\min\{\kappa_1,\kappa_2\}$ proves both decay estimates.
\end{proof}

\subsection{Proof of the periodic theorem}\label{sec:completion}
We combine the finite-time continuation bounds with the stability estimate from the entry time onward.
\begin{proof}[Proof of Theorem~\ref{thm:periodic}]
The global solution from Section~\ref{sec:continuation} satisfies
\eqref{eq:strong-regularity} on each finite interval. By Lemma~\ref{lem:good-times}
with $\varepsilon=\varepsilon_*$, there is $t_*<L_*:=L_{\varepsilon_*}$ at
which Proposition~\ref{prop:tail-stability} applies. With horizon $L_*+1$,
Proposition~\ref{prop:upper-poly} for $1<\gamma\le3$ and
Proposition~\ref{prop:upper-iso} for $\gamma=1$ give the upper density
bound on $[0,L_*]$. Proposition~\ref{prop:lower} then gives the lower bound
there. Write $[c_*,R_*]$ for the resulting density interval. By
\eqref{eq:tail-density}, $[M/2,3M/2]$ suffices on $[t_*,\infty)$. Thus
$\min\{c_*,M/2\}\le\rho\le\max\{R_*,3M/2\}$ for all time.

For $t\ge t_*$, \eqref{eq:tail-log}, $u=v-\nabla\log\rho$, and the mean-zero Poisson estimate give
\begin{equation}\label{eq:decay-final}
 \|\rho-M\|_{H^3}+\|u-U\|_{H^2}+\|V\|_{H^5}
 \le C\bigl(\|\rho-M\|_{H^3}+\|v-U\|_{H^2}\bigr)
 \le C A_*^{1/2} e^{-\kappa(t-t_*)/2}.
\end{equation}
On the unit torus $c_0\le M\le C_0$ and $|U|\le\|u_0\|_\infty$.
The initial bounds \eqref{eq:initial-data}, the formula
$v_0=u_0+\nabla\rho_0/\rho_0$ and the periodic Poisson estimate of
Lemma~\ref{lem:poisson-estimates} bound $E_u(0)+E_v(0)$ in terms of the data.
The constants in Lemma~\ref{lem:good-times}, Proposition~\ref{prop:tail-stability},
Propositions~\ref{prop:upper-poly} and~\ref{prop:lower}, and
Proposition~\ref{prop:global} can be chosen uniformly for
$M\in[c_0,C_0]$, as their proofs use only upper bounds for $M$ and $M^{-1}$
at fixed $\gamma,\ell$. Thus \eqref{eq:tail-entry-time} gives a bound for
$L_*$ in terms of the data, and
Proposition~\ref{prop:global} on $[0,L_*]$ bounds $A_*$ by the same data.
Choose $T_0$ to be the resulting data-dependent upper bound for $L_*$.
Since $t_*<L_*\le T_0$, \eqref{eq:decay-final} gives
\eqref{eq:main-exponential-decay} for $t\ge T_0$ with $\sigma=\kappa/2$
and a prefactor controlled by the data. The rate can
be chosen in terms of $M,\gamma,\ell$ and the torus.
\end{proof}

\section{Time-uniform density bounds on the whole space}\label{sec:main-proof}

Proposition~\ref{prop:global} gives global existence and uniqueness.
To obtain \eqref{eq:uniform}, we use the integrated dissipation to
control the data for the finite-interval estimates at a time of bounded
distance from each observation time.

\begin{proof}[Proof of Theorem~\ref{thm:cauchy}]
Fix an observation time $t$. When $t$ is sufficiently large, we first
find an earlier time $t_0$ at which $\rho$ is close to $b$ in
$L^\infty$ and $m$ is bounded in $L^6$. Set
\begin{equation}\label{eq:window-parameters}
 \varepsilon_b=\min\Bigl\{1,\frac{\sqrt b}{4S_2(1+b^{-1/2})},
                   \frac1{2S_1\sqrt{2b}}\Bigr\},\qquad
 T_{\rm w}=1+\max\Bigl\{1,\frac{(7+\ell^2)\budget}{\varepsilon_b^2}\Bigr\}.
\end{equation}
Let $D(\tau)$ denote the integrand in \eqref{eq:dissipation-budget}.
It is continuous on finite intervals by \eqref{eq:regularity},
\eqref{eq:compact-positive} and \eqref{eq:effective-variables}. For any
$t\ge T_{\rm w}$ there is a $t_0\in[t-T_{\rm w},t-1]$ such that
\begin{equation}\label{eq:small-diss-time}
 D(t_0)\le\frac1{T_{\rm w}-1}\int_{t-T_{\rm w}}^{t-1}D(\tau)\dd\tau
 \le\frac{(7+\ell^2)\budget}{T_{\rm w}-1}\le\varepsilon_b^2.
\end{equation}
At this time, $f=\sqrt\rho-\sqrt b=r/(\sqrt\rho+\sqrt b)$ satisfies
$\|\Delta f\|_2=\|\nabla^2s\|_2$ and, by the elliptic estimate
$\|f\|_{H^2}\le C(\|f\|_2+\|\Delta f\|_2)$, which follows frequency by
frequency from $(1+|\xi|^2)^2\le C(1+|\xi|^4)$; the constant $C$ is absorbed
into $S_2$ below. Since
$\|f\|_2\le b^{-1/2}\|r\|_2$, we obtain
\[
 \|f\|_\infty\le S_2\bigl(b^{-1/2}\|r\|_2+\|\nabla^2s\|_2\bigr)
 \le\sqrt b/4.
\]
Consequently $9b/16\le\rho(t_0)\le25b/16$. Also, since $v$ vanishes at
infinity and $\rho(t_0)\le25b/16<2b$,
\[
 \|m(t_0)\|_6\le 2b S_1\|\nabla v(t_0)\|_2
 \le\tfrac83 S_1\sqrt b\,\|\sqrt{\rho(t_0)}\,\nabla v(t_0)\|_2\le1.
\]

We propagate these estimates long enough for the momentum equation to
give an $L^\infty$ bound. On any interval starting at $t_0$ on which
$b/4\le\rho\le3b$ and $\|m\|_6\le2$, we have $\|r\|_\infty\le2b$.
Recall that $K_r$ denotes the right-hand side of \eqref{eq:charge-budget}. Since
$|P(\rho)-P(b)|\le\gamma(3b)^{\gamma-1}|r|$ and
$\|r\|_6\le\|r\|_2^{1/3}\|r\|_\infty^{2/3}$, while
$\|E\|_6\le S_1\ell^{-2}\|r\|_2$ and \eqref{eq:newton-field} bounds
$\|E\|_\infty$, we have
\begin{align*}
 \|\Pi\|_6+\|\Pi\|_\infty
 &\le\gamma(3b)^{\gamma-1}\bigl(K_r^{1/3}(2b)^{2/3}+2b\bigr),\\
 \|\rho E\|_6+\|\rho E\|_\infty
 &\le3b\ell^{-2}\bigl[(S_1+(4\pi)^{-1/2})K_r+2b\bigr].
\end{align*}
Choose $\bar B\ge1$ bounding both right-hand sides. Fix
$0<\delta\le1/2$ such that
\begin{equation}\label{eq:restart-delta}
 \frac{16C_{\rm heat}}b\delta^{1/4}\le\tfrac14,\qquad
 C_{\rm heat}\bar B\sqrt\delta\le\tfrac18.
\end{equation}
These choices are independent of $t$ and $t_0$.

For $0<h\le\delta$, while the stated bounds hold, the continuity equation in
\eqref{eq:heat-mass} and \eqref{eq:heat-kernel} give
\[
\rho(t_0+h)=G_h*\rho(t_0)
-\int_0^h\sum_j\partial_jG_{h-\tau}*m_j(t_0+\tau)\dd\tau ,
\]
and since $\|\partial_jG_{h-\tau}*m_j\|_\infty
\le(C_{\rm heat}/8)(h-\tau)^{-3/4}\|m\|_6$, the gradient heat estimate from
$L^6$ to $L^\infty$ gives
\begin{equation}\label{eq:window-density}
 \|\rho(t_0+h)-G_h*\rho(t_0)\|_\infty
 \le \frac{C_{\rm heat}}8\int_0^h 2(h-\tau)^{-3/4}\dd\tau
 =C_{\rm heat}h^{1/4}\le b/64.
\end{equation}
The heat average $G_h*\rho(t_0)$ is a convex combination of values of
$\rho(t_0)$, hence it lies between $9b/16$ and $25b/16$. Since
$\|m\otimes m/\rho\|_3\le(4/b)\|m\|_6^2\le16/b$, the momentum equation,
with its differentiated kernels acting from $L^3$ to $L^6$ on the convective
flux and from $L^6$ to $L^6$ on the pressure, yields
\begin{equation}\label{eq:window-six}
 \|m(t_0+h)\|_6\le1+\frac{8C_{\rm heat}}b h^{1/4}
            +\tfrac14 C_{\rm heat}\bar B\sqrt h+\bar Bh<\tfrac32.
\end{equation}
In the last inequality we used \eqref{eq:restart-delta},
$C_{\rm heat}\ge8$ and $h\le1$. Both estimates strictly improve the
bounds $b/4\le\rho\le3b$, $\|m\|_6\le2$. These bounds are strict at
$t_0$, and $r$ and $m$ are continuous in $L^\infty$ and $L^6$,
respectively. A first-exit argument now gives both bounds on
$[t_0,t_0+\delta]$.

For this particular interval the quantity
$X=\sup_{0<h\le\delta}h^{1/4}\|m(t_0+h)\|_\infty$ is finite by
$H^2$ continuity. To bound it independently of $t_0$, use
$\|m\otimes m/\rho\|_6\le (8/b)\|m\|_\infty$ and apply
\eqref{eq:heat-integral} to the momentum. The $L^6$ initial trace
contributes at most $\|G_1\|_{6/5}$ after multiplication by $h^{1/4}$.
Thus
\begin{align*}
 h^{1/4}\|m(t_0+h)\|_\infty
 &\le\|G_1\|_{6/5}
 +\frac{C_{\rm heat}}b h^{1/4}X
   \int_0^h(h-\tau)^{-3/4}\tau^{-1/4}\dd\tau\\
 &\quad+h^{1/4}\bar B\bigl(\tfrac14 C_{\rm heat}\sqrt h+h\bigr).
\end{align*}
The change of variables $\tau=hy$ reduces the integral to
$\int_0^1(1-y)^{-3/4}y^{-1/4}\dd y$. Splitting at $1/2$ bounds it by
$4/3+4<8$. Taking the supremum and using \eqref{eq:restart-delta}
gives $X\le C_{\rm heat}+\frac18X+\frac14$, hence $X\le2C_{\rm heat}$.
At $a=t_0+\delta$ we have proved
\begin{equation}\label{eq:restart-trace}
 b/4\le\rho(a)\le3b,\qquad
 \|v(a)\|_\infty\le\frac{8C_{\rm heat}}{b\delta^{1/4}}.
\end{equation}

Given $t>0$, take $a=0$ if $t\le T_{\rm w}$, and take $a=t_0+\delta$
as above otherwise. In the latter case
$1-\delta\le t-a\le T_{\rm w}-\delta$. In both cases $0<t-a\le T_{\rm w}$,
and the data bounds at $a$ can be taken to be
\begin{align*}
 c_{\rm in}&=\min\{c_0,b/4\},\qquad R_{\rm in}=\max\{C_0,3b\},\\
 W_{\rm in}&=\max\left\{S_2\bigl(\|u_0\|_{H^2}
                 +c_0^{-1}\|r_0\|_{H^3}\bigr),
                   \frac{8C_{\rm heat}}{b\delta^{1/4}}\right\}.
\end{align*}
Proposition~\ref{prop:upper} on $[a,t]$ gives the same
$M_{T_{\rm w}}$ for every $t$. \eqref{eq:field-finite-bound}
and Lemma~\ref{lem:large-moments} then give a common $A_{T_{\rm w}}$,
so Proposition~\ref{prop:lower} gives a common $c_{T_{\rm w}}>0$.
This covers every $t\ge T_{\rm w}$; for $0\le t\le T_{\rm w}$ the same
bounds follow from the trace at $a=0$, and the case $t=0$ is the initial
datum. We may therefore take
$C_{\rm den}=\max\{1,M_{T_{\rm w}},c_{T_{\rm w}}^{-1}\}$ in
\eqref{eq:uniform}. Its dependence on the data is the one stated in
Theorem~\ref{thm:cauchy}.
\end{proof}

\appendix
\section{Local existence and uniqueness}\label{app:local}

The density perturbation, effective momentum and electric field form a closed
semilinear heat system. Its local existence time depends on an $H^2$
bound for the three variables and a positive lower bound for the density.
The density gains the additional derivative required to recover the
physical velocity from the effective momentum when its initial trace
belongs to $H^3$.

\begin{proposition}\label{prop:local-existence}
Let $K\ge0$, $b,\ell>0$, $0<c\le b$ and $1\le\gamma\le3$; the condition
$c\le b$ makes the constant state $(b,0)$ admissible. There is an
$h=h(K,c,b,\ell,\gamma)>0$ with the following property. Suppose $a\ge0$
and $\mathcal U_a=(r_a,m_a,E_a)\in H^2$ satisfies
\[
 \|\mathcal U_a\|_{H^2}\le K,\qquad b+r_a\ge c,\qquad
 \nabla\times E_a=0,\qquad \Div E_a=-\ell^{-2}r_a.
\]
On $\T$ assume in addition that $r_a$ and $E_a$ have zero mean. On $J=[a,a+h]$, equations \eqref{eq:heat-mass}--\eqref{eq:heat-momentum} and
\eqref{eq:field-heat} have a unique mild solution $\mathcal U\in C(J;H^2)$
with trace $\mathcal U_a$, in the class of mild solutions for which $b+r$ is
positive at every point of $\Omega\times J$. It satisfies
$\rho=b+r\ge c/2$, $\nabla\times E=0$ and $\Div E=-\ell^{-2}r$.

If in addition $r_a\in H^3$, then $\rho=b+r$, $u=(m-\nabla r)/\rho$ and
the normalized potential $V$ form the unique positive strong solution on $J$
in the sense of Definition~\ref{def:positive-strong-solution}, with density
$b+r_a$ and velocity $(m_a-\nabla r_a)/(b+r_a)$ at time $a$.
The solution has the regularity
\eqref{eq:strong-regularity} on $\T$ and \eqref{eq:regularity} on $\R^3$, with time interval $J$.
\end{proposition}

\begin{proof}
Let $\eta=\min\{1,c/(2S_2)\}$. In $C(J;H^2)$ take the closed
ball $\mathcal B$ of radius $\eta$ about the free heat evolution
$t\mapsto e^{(t-a)\Delta}\mathcal U_a$.
The free evolution has norm at most $K$, and its density satisfies
$b+e^{(t-a)\Delta}r_a\ge c$ by positivity of the heat kernel.
Consequently every $W\in\mathcal B$ has norm at most $K+1$ and density
at least $c/2$, regardless of $h$: indeed
$\|W\|_{H^2}\le K+\eta\le K+1$, and, since $H^2\hookrightarrow L^\infty$
with constant $S_2$ and $\eta\le c/(2S_2)$,
$b+W_\rho\ge b+e^{(t-a)\Delta}r_a-S_2\eta\ge c/2$.

We record the nonlinear estimates for the contraction. The product
rule and the embeddings $H^2\hookrightarrow L^\infty\cap W^{1,4}$ imply
\begin{equation}\label{eq:local-product}
 \|gw\|_{H^2}\le C\|g\|_{H^2}\|w\|_{H^2}.
\end{equation}
Let $F$ be a pointwise component of the maps in \eqref{eq:heat-maps}.
It vanishes at zero, and its first three derivatives are bounded on
$\{|\mathcal U|\le S_2(K+1),\, b+r\ge c/2\}$. This set is convex and contains
the origin. For $W,\widetilde W\in\mathcal B$, write
$W_\theta=\widetilde W+\theta(W-\widetilde W)$. The mean value formula
and the chain rule give, at each time,
\begin{align*}
 \|DF(W_\theta)-DF(0)\|_{H^2}
 &\le C\bigl(\|W_\theta\|_{H^2}+\|\nabla W_\theta\|_4^2\bigr),\\
 F(W)-F(\widetilde W)
 &=DF(0)(W-\widetilde W)\\
 &\quad+\int_0^1\!\bigl(DF(W_\theta)-DF(0)\bigr)(W-\widetilde W)\dd \theta.
\end{align*}
By \eqref{eq:local-product}, the resulting composition estimate is
\begin{equation}\label{eq:local-composition}
 \|F(W)-F(\widetilde W)\|_{H^2}\le C\|W-\widetilde W\|_{H^2},
 \qquad \|F(W)\|_{H^2}\le C\|W\|_{H^2}.
\end{equation}
For the second estimate use the segment from $0$ to $W$.
The remaining component is linear, with
$\|\ell^{-2}Qm\|_{H^2}\le\ell^{-2}\|m\|_{H^2}$. Thus a single constant
$B_{\rm loc}\ge1$, depending only on $K,c,b,\ell,\gamma$, bounds the
norms and the Lipschitz constants of the flux and source on $\mathcal B$.

Define $\mathcal T$ by the right-hand side of \eqref{eq:heat-integral}
with initial trace $\mathcal U_a$, flux $\mathcal F(W)$ and source $\mathcal S(W)$.
Since the Fourier multipliers satisfy $e^{-d|\xi|^2}\le1$ and
$\sup_\xi|\xi|e^{-d|\xi|^2}\le d^{-1/2}$ for $d>0$, the convolution
operators obey $\|G_d*f\|_{H^2}\le C\|f\|_{H^2}$ and
$\|\nabla G_d*f\|_{H^2}\le Cd^{-1/2}\|f\|_{H^2}$ for $0<d\le1$; hence the
last
$\varepsilon$ units of its time integral have $H^2$ norm at most
$B_{\rm loc}\int_0^\varepsilon(1+s^{-1/2})\dd s
=B_{\rm loc}(2\sqrt\varepsilon+\varepsilon)$. Away from the upper
endpoint, strong continuity of the heat semigroup gives continuity of
the integral in $H^2$. In particular $\mathcal T W\in C(J;H^2)$ and
\begin{align*}
 \|\mathcal T W-e^{(\,\cdot-a)\Delta}\mathcal U_a\|_{C(J;H^2)}
   &\le B_{\rm loc}(2\sqrt h+h),\\
 \|\mathcal T W-\mathcal T\widetilde W\|_{C(J;H^2)}
   &\le B_{\rm loc}(2\sqrt h+h)\|W-\widetilde W\|_{C(J;H^2)}.
\end{align*}
Choose $h\in(0,1]$ so that
\begin{equation}\label{eq:local-lifespan}
 B_{\rm loc}(2\sqrt h+h)\le\min\{\eta,1/2\}.
\end{equation}
Banach's fixed point theorem gives a mild solution in $\mathcal B$.
On $\T$, the heat equations preserve the zero means of $r$ and $E$. The homogeneous heat equations \eqref{eq:constraint-heat} propagate
the two field constraints. Any other positive-density $C(J;H^2)$ mild
solution with the same trace belongs to the corresponding ball on a
short initial interval by continuity. At each common trace the density
has a positive infimum by compactness on the torus or convergence to $b$ at infinity. Restarting the contraction at the last time of agreement, if
it were earlier than the right endpoint of $J$, proves uniqueness on $J$.

The flux and source of the fixed point belong to $C(J;H^2)$. The
linear estimate \eqref{eq:linear-integrated} at $k=2$ controls
$m$ in $L^2(J;H^3)$. If $r_a\in H^3$, the same estimate at $k=3$ for
$r_t-\Delta r=-\Div m$ gives
\begin{equation}\label{eq:local-spaces}
 \mathcal U\in C(J;H^2)\cap L^2(J;H^3),\qquad
 r\in C(J;H^3)\cap L^2(J;H^4).
\end{equation}
The norm $\|r_a\|_{H^3}$ enters this estimate only after the lifespan
has been fixed in \eqref{eq:local-lifespan}. The equations then give
\[
 r_t\in C(J;H^1)\cap L^2(J;H^2),\qquad
 m_t\in C(J;L^2)\cap L^2(J;H^1).
\]
Since $\rho\ge c/2$ and $r\in C(J;H^3)$, the chain rule yields
$1/\rho-1/b\in C(J;H^3)$. Indeed, $H^3\hookrightarrow W^{1,\infty}$
controls the factors $\nabla r$ in the mixed and cubic terms of the
third derivative. Multiplication by $1/\rho$ is therefore bounded on
$H^2$ and $H^3$, uniformly on $J$. The formulas $u=(m-\nabla r)/\rho$
and $v=m/\rho$ give
$u,v\in C(J;H^2)\cap L^2(J;H^3)$. Moreover, in distributions,
\begin{equation}\label{eq:velocity-derivs}
 \rho u_t=m_t-\nabla r_t-ur_t,\qquad \rho v_t=m_t-vr_t.
\end{equation}
For $w\in H^2$, $z\in H^1$, the embeddings into $L^\infty$ and $L^6$,
respectively, and $\|\nabla w\|_3\le C\|w\|_{H^2}$ give
$\|wz\|_{H^1}\le C\|w\|_{H^2}\|z\|_{H^1}$. Thus $ur_t,vr_t\in C(J;H^1)$.
It follows from \eqref{eq:velocity-derivs} and \eqref{eq:local-spaces} that
\begin{equation}\label{eq:velocity-traces}
 u_t,v_t\in C(J;L^2)\cap L^2(J;H^1).
\end{equation}

For the electric field, the Fourier identity \eqref{eq:gauss-symbol}
controls high frequencies by the density, while the mild solution
already controls the low frequencies in $L^2$. Splitting at $|\xi|=1$
gives, for every integer $k\ge0$ with $r\in H^k$,
\begin{equation}\label{eq:field-smoothing}
 \|E\|_{H^{k+1}}^2\le 2^{k+1}\|E\|_2^2+2\ell^{-4}\|r\|_{H^k}^2.
\end{equation}
Apply \eqref{eq:field-smoothing} at $k=3$ to differences of time
slices and at $k=4$ under the time integral. By
\eqref{eq:local-spaces}, $E\in C(J;H^4)\cap L^2(J;H^5)$. On the torus the zero spatial mean of $V$ then gives $V\in C(J;H^5)\cap L^2(J;H^6)$.

The identities \eqref{eq:effective-flux} and
\eqref{eq:heat-mass}--\eqref{eq:field-heat}, read in reverse, give
\eqref{eq:mass}--\eqref{eq:poisson}. The formulas for $\rho,u$ give the
specified traces. In the whole space, $r,u\in H^2\hookrightarrow C_0(\R^3)$ at every
time, so they have the prescribed limits at infinity.
In $\R^3$, \eqref{eq:potential-normalized} supplies the potential vanishing
at infinity. Conversely, any positive strong solution on $J$ in the sense
of Definition~\ref{def:positive-strong-solution} with these data satisfies
$m=\rho u+\nabla r\in C(J;H^2)$ by \eqref{eq:local-product}. The
distributional Poisson equation gives the
field constraints at each time by continuity of $r$ and $E$. Applying
\eqref{eq:field-smoothing} to differences of time slices then gives
$E\in C(J;H^4)$. Thus $\mathcal U\in C(J;H^2)$. The compactness argument
following Definition~\ref{def:positive-strong-solution} gives a positive lower
bound for $\rho$ on $J$.
The computations leading to \eqref{eq:heat-mass}--\eqref{eq:heat-momentum}
and \eqref{eq:field-heat} give the heat
system in distributions. Duhamel's formula \eqref{eq:heat-integral}
and uniqueness of the mild solution now determine
$\rho$ and $u$. The normalization determines $V$.

Finally, for the data of Theorems~\ref{thm:periodic} and~\ref{thm:cauchy}, the product estimate
\eqref{eq:local-product} applied to $m_0=bu_0+r_0u_0+\nabla r_0$ and
\eqref{eq:field-smoothing} at $k=1$ give
\begin{equation}\label{eq:initial-effective-data}
 \begin{aligned}
 \|m_0\|_{H^2}
 &\le (b+C\|r_0\|_{H^2})\|u_0\|_{H^2}+\|r_0\|_{H^3},\\
 \|E_0\|_{H^2}^2&\le4\|E_0\|_2^2+2\ell^{-4}\|r_0\|_{H^1}^2.
 \end{aligned}
\end{equation}
On the torus, the spectral gap gives $\|E_0\|_2\le (2\pi\ell^2)^{-1}\|r_0\|_2$. In the whole space \eqref{eq:data-field} gives the corresponding low-frequency bound. Thus both initial states have the required $H^2$ regularity.
\end{proof}

\section*{Acknowledgments}
	X. Huang is partially supported by Chinese Academy of Sciences Project for Young Scientists in Basic Research (Grant No. YSBR-031), National Natural Science Foundation of China (Grant Nos. 12494542, 11688101) and National Key R\&D Program of China (Grant No. 2021YFA1000801). 
	
	\smallskip
	\noindent\textbf{AI assistance.} The authors acknowledge the use of AI to explore proof
	strategies, check calculations, and polish the manuscript. In particular, AI suggested the
	mixed density--velocity moment used to obtain the additional density integrability at the
	endpoint $\gamma=3$ and helped verify the exponential decay on the torus. The authors take
	full responsibility for the mathematical arguments and the final manuscript.
	
	\smallskip
	\noindent\textbf{Data availability statement.} Data sharing is not applicable to this article.
	
	\smallskip
	\noindent\textbf{Conflict of interest.} The authors declare that they have no conflict of interest.

\end{document}